\documentclass[a4paper,12pt]{article}
\usepackage[T2A]{fontenc}
\usepackage[cp1251]{inputenc}
\begin{document}

\author{S.V. Ludkowski}

\title{Octonion random functions and integration of stochastic PDEs.}

\date{29 September 2020}
\maketitle

\begin{abstract}
In the article random functions in modules over the octonion algebra
and Cayley-Dickson algebras are investigated. For their study
transition measures with values in the octonion algebra and
Cayley-Dickson algebras are used. Stochastic integrals over these
algebras are studied. They are applied to integration of stochastic
PDEs. This approach permits subsequently to analyze and integrate
PDEs of orders higher than two of different types including
parabolic, elliptic and hyperbolic.

\footnote{key words and phrases: Octonion
Algebra; Cayley-Dickson Algebra; Random Function; Integration; Stochastic PDE 
address: Dep. Appl. Mathematics,  MIREA - Russian Technological
University, av. Vernadsky 78, Moscow, 119454, Russia \par e-mail:
sludkowski@mail.ru}
\end{abstract}

\section{Introduction.}
Over the complex field Feynman integrals and quasi-measures appeared
to be very important in mathematics, mathematical physics, quantum
mechanics, quantum field theory and partial differential equations
(PDEs) (see, for example,
\cite{dalfom,gulcastb,johnlap,kimfi09,nicolam16}). It is worthwhile
to note that in stochastic analysis measures with values in matrix
algebras or operator algebras on Hilbert spaces are frequently
studied and used for solution or analysis of PDEs
\cite{dalfom,gihmskorb,gulcastb,meashand,johnlap}.
\par But there are restrictions for these approaches,
because the Feynman integral works for partial differential
operators (PDOs) of order not higher than two. Indeed, it is based
on complex modifications of Gaussian measures. Nevertheless, if a
characteristic function $\phi (t)$ of a measure has the form $\phi
(t)=\exp (Q(t))$, where $Q(t)$ is a polynomial, then its degree is
not higher than two according to the Marcinkievich theorem (Ch. II,
\S 12 in \cite{shirb11}). \par On the other side, hypercomplex
numbers open new opportunities in these areas. For example, Dirac
used the complexified quaternion algebra ${\bf H}_C$ for a solution
the Klein-Gordon hyperbolic PDE of the second order with constant
coefficients \cite{dirac}. This is important in spin quantum
mechanics. It was proved in \cite{ludrend14} that in many variants,
it is possible to reduce a PDE problem to a subsequent solution of
PDEs of order not higher than two with Cayley-Dickson coefficients.
In general the complex field is insufficient for this purpose.
\par On the other hand, the Cayley-Dickson algebras ${\cal A}_r$ over
the real field $\bf R$ are natural generalizations of the complex
field, where ${\cal A}_2 = \bf H$ is the quaternion skew field,
${\cal A}_3=\bf O$ denotes the octonion algebra, ${\cal A}_0=\bf R$,
${\cal A}_1=\bf C$. Then each subsequent algebra ${\cal A}_{r+1}$ is
obtained from the preceding algebra ${\cal A}_r$ by the doubling
procedure using the doubling generator \cite{baez,dickson,kansol}.
\par They are widely applied in PDEs, non-commutative analysis,
mathematical physics, quantum field theory, hydrodynamics,
industrial and computational mathematics, non-commutative geometry
\cite{dirac,emch,guespr,guesprqa,guetze,lawmich,luoyst,ludjms7,lufsqv,
lufoclg}.
\par Previously measures with values in the complexified
Cayley-Dickson algebra ${\cal A}_{r,C}$ were studied in
\cite{luocmeasaaca20}. They appear naturally while a solution of the
second order hyperbolic PDE with Cayley-Dickson coefficients. In
this work the results and notation of \cite{luocmeasaaca20} are
used. This article is devoted to a realization of the plan
formulated in the preceding cited work: for solution and analysis of
PDEs of orders higher than two to extend Feynman integrals and
quasi-measures from spaces over the complex field onto modules over
the Cayley-Dickson algebras.
\par In this paper random functions in modules
over the complexified octonion algebra ${\bf O}_C={\cal A}_{3,C}$
and the complexified Cayley-Dickson algebras ${\cal A}_{r,C}$ are
investigated. For their study transition measures with values in the
complexified octonion algebra and the complexified Cayley-Dickson
algebras are used. An existence of random functions and Markov
processes in modules over the complexified Cayley-Dickson algebra
${\cal A}_{r,C}$ is studied in Theorem 2.8, Corollary 2.9.
Stochastic integrals of such random functions and acting on them
operators are investigated in Theorems 2.14, 2.15, 2.17, 2.18.
Properties of stochastic integrals over ${\cal A}_{r,C}$ are
described in Propositions 2.20-2.22. In Theorem 2.27 their
stochastic continuity is studied. Necessary specific definitions are
given. Notation is described in remarks. Lemmas 2.12, 2.13, 2.25,
2.26 are given in order to prove the theorems and propositions.
These lemmas concern estimates of stochastic integrals over ${\cal
A}_{r,C}$. In Theorems 2.29, 2.31 and Corollary 2.30 solutions of
stochastic PDEs with random functions and operators in modules over
${\cal A}_{r,C}$ are scrutinized.
\par A formula number $(n)$ in the same subsection m is referred as
$(n)$, in another subsection as m$(n)$.
\par Main results of this work are
obtained for the first time. The obtained results open new
opportunities for subsequent studies of PDEs and their solutions
including that of hyperbolic type and parabolic type with hyperbolic
and elliptic terms of orders two or higher, related to them random
functions.

\section{Octonion random functions.}

\par {\bf Definition 2.1.} Suppose that $\bf \Lambda $ is an additive
group contained in $\bf R$. Suppose also that $T$ is a subset in
$\bf \Lambda $ and containing a point $t_0$. Let $X_t=X$ be a
locally $\bf R$-convex space which be also a two-sided ${\cal
A}_{r,C}$-module for each $t\in T$, where $2\le r<\infty $. Put
$$(\tilde{X}_T, {\tilde {\sf U}}) := \prod_{t\in T}(X_t,{\sf U}_t)$$
for the product of measurable spaces, where ${\sf U}_t$ is the Borel
$\sigma $-algebra of $X_t$, $ ~ \tilde {\sf U}$ is an algebra of
cylindrical subsets of $\tilde{X}_T$ generated by projections
$\tilde \pi _q: \tilde{X}_T\to X^q$, where
$X^q:=\mbox{}_l\prod_{t\in q}X_t$ is a left ordered direct product,
$q\subset T$ is a finite subset of $T$, $~X^{ \{ t \} } =X_t$,
$~X^{t_1,...,t_{n+1}} = X_{t_{n+1}}\times (X^{t_1,...,t_n})$ for
each $t_1<...<t_{n+1}$ in $T$.

A function $P(t_1,x_1,t_2,A)$ with values in the complexified
Cayley-Dickson algebra ${\cal A}_{r,C}$ for each $t_1< t_2\in T$,
$x_1\in X_{t_1}$, $A\in {\sf U}_{t_2}$ is called a transition
measure if it satisfies the following conditions:

$(1)\mbox{ the set function }\nu _{x_1,t_1,t_2}(A)
:=P(t_1,x_1,t_2,A)\mbox{ is a } \mbox{measure on }$ $(X_{t_2},{\sf
U}_{t_2});$

$(2)\mbox{ the function }\alpha _{t_1,t_2,A}(x_1)
:=P(t_1,x_1,t_2,A)\mbox{ of the variable } x_1$ $\mbox{ is  } {\sf
U}_{t_1}$-mea\-su\-rable, that is, $\alpha _{t_1,t_2,A}^{-1}({\cal
B}({\cal A}_{r,C}))\subset {\sf U}_{t_1}$;
$$(3)\quad P(t_1,x_1,t_2,A)=\int_{X_z}
P(t,y,t_2,A)P(t_1,x_1,t,dy)\mbox{ for each }t_1<t<t_2\in T$$ so that
$P(t,y,t_2,A)$ as the function by $y$ is in $L^1((X_t,{\sf U}_t),\nu
_{x_1,t_1,t},{\cal A}_{r,C})$. A transition measure
$P(t_1,x_1,t_2,A)$ is called unital if
$$(4)\quad P(t_1,x_1,t_2,X_{t_2})=1\mbox{ for each }t_1<t_2\in T.$$

\par Then for each finite set $q=(t_0,t_1,\dots,t_{n+1})$
of points in $T$ such that $t_0<t_1<...<t_{n+1}$ there is defined a
measure in $X^g$
$$(5)\quad \mu ^q_{x_0}(D) = \int_D
\mbox{ }_l\prod_{k=1}^{n+1}P(t_{k-1},x_{k-1}, t_k,dx_k), \mbox{
}D\in {\sf U}^g := \mbox{}_l\prod_{t\in g}{\sf U}_t,$$ where
$g=q\setminus \{ t_0 \}$, variables $x_1,\dots ,x_{n+1}$ are such
that $(x_1,\dots ,x_{n+1})\in D$, $ ~ x_0\in X_{t_0}$ is fixed.
\par Let the transition measure $P(t,x_1,t_2,dx_2)$ be
unital. Then for the product $D=D_2\times (X_{t_j}\times D_1)$,
where $D_1\in \mbox{}_l\prod_{i=1}^{j-1}{\sf U}_{t_i}$, $ ~ D_2\in
\mbox{}_l\prod_{i=j+1}^{n+1} {\sf U}_{t_i}$,  the equality
$$(6)\quad \mu ^q_{x_0}(D) = \int_{D_2\times D_1}\bigg[
\mbox{}_l\prod_{k=j+1}^{n+1}P(t_{k-1},x_{k-1},t_k,dx_k)\bigg]$$
$$ \times \bigg[ \bigg[\int_{X_{t_j}}P(t_{j-1},x_{j-1},t_j,dx_j)
\bigg[\mbox{}_l\prod_{k=1}^{j-1}P( t_{k-1},x_{k-1},t_k,dx_k) \bigg]
\bigg]
 = \mu ^r_{x_0}(D_2\times D_1)$$
is fulfilled, where $r=q\setminus \{ t_j \} .$ Equation $(6)$
implies that
$$(7)\quad [\mu ^q_{x_0}]^{\pi ^q_v}=\mu ^v_{x_0}$$
for each $v<q$, where finite sets are ordered by inclusion: $v<q$ if
and only if $v\subset q$, where $\pi ^q_w: X^g\to X^w$ is the
natural projection, $g=q\setminus \{ t_0 \} ,$ $ ~ w=v\setminus \{
t_0 \} .$
\par Denote by $\Upsilon _T$ the family of all finite linearly ordered
subsets $q$ in $T$ such that $t_0\in q\subset T$, $v\le q\in
\Upsilon _T$, $\pi _q: \tilde X_T\to X^g$ is the natural projection,
$g=q\setminus \{ t_0 \} $. Hence Conditions $(4)$, $(5)$, $(7)$
imply that: $\{ \mu ^q_{x_0}; \pi ^q_v; \Upsilon _T \} $ is the
consistent family of measures. It induces a cylindrical distribution
$\tilde \mu _{x_0}$ on the measurable space $(\tilde X_T, \tilde
{\sf U})$ such that
$$(8)\quad \tilde \mu _{x_0}(\pi _q^{-1}(D))=\mu ^q_{x_0}(D)$$  for each
$D\in {\sf U}^g$.
\par The cylindrical distribution given by Formulas $(1)$-$(5)$, $(7)$,
$(8)$ is called the ${\cal A}_{r,C}$-valued Markov distribution with
time $t$ in $T$.

\par {\bf Remark 2.2.} Let $X_t=X$ for each $t\in T$, ${\tilde X}_{t_0,x_0}:=
\{ x\in {\tilde X}_T: $ $x(t_0)=x_0 \} .$ Put $\bar{\pi }_q:$
$x\mapsto x_q$ for each $x=x(t)$ in ${\tilde X}_T$, where $x_q$ is
defined on $q=(t_0,\dots ,t_{n+1})\in \Upsilon _T$ such that
$x_q(t)=x(t)$ for each $t\in q$. To an arbitrary function $F:
{\tilde X}_T\to {\cal A}_{r,C}^l$ a function can be posed
$(S_qF)(x):=F(x_q)=F_q(y_0,\dots ,y_n),$ where $y_j=x(t_j)$, $F_q:
X^q\to {\cal A}_{r,C}^l$, $~l\in \bf N$. Put
$${\sf F}:= \{ F| F: {\tilde X}_T\to {\cal A}_{r,C}^l, S_qF\mbox{ is }{\sf
U}^q-\mbox{measurable for each } q\in \Upsilon _T \} .$$ If $F\in
\sf F$, $ ~ \tau =t_0\in q$, $ ~ t_0<t_1<...<t_{n+1}$, then the
integral
$$(1)\quad J_q(F)=\int_{X^q}(S_qF)
(x_0,\dots ,x_n)\mbox{ }_l\prod_{k=1}^{n+1}P(t_{k-1},
x_{k-1},t_k,dx_k)$$ can be defined whenever it converges.

\par {\bf Definition 2.3.} A function $F$ is called
integrable relative to a Markov cylindrical distribution $\mu
_{x_0}$ if the limit
$$(1)\quad \lim_{q\in \Upsilon _T} J_q(F) =: J(F)$$ along the generalized net by
finite subsets $q=(t_0,\dots ,t_{n+1})\in \Upsilon _T$ of $T$
exists. This limit is called a functional integral relative to the
Markov cylindrical distribution:
$$(2)\quad J(F)=\int_{{\tilde X}_{t_0,x_0}}F(x)\mu _{x_0}(dx).$$

\par {\bf Remark 2.4. Spatially
homogeneous transition measure.} Suppose that $P(t,A)$ is an ${\cal
A}_{r,C}$-valued measure on $(X,{\sf U})$ for each $t\in T$ such
that $A-x\in \sf U$ for each $A\in \sf U$ and $x\in X$, where $A\in
\sf U$, $X$ is a locally $\bf R$-convex space which is also a
two-sided ${\cal A}_{r,C}$-module, $\sf U$ is an algebra of subsets
of $X$. Suppose also that $P$ is a spatially homogeneous transition
measure:
$$(1)\quad P(t_1,x_1,t_2,A)=P(t_2-t_1,A-x_1)$$
for each $A\in \sf U$, $t_1<t_2 \in T$ and $t_2-t_1\in T$ and every
$x_1\in X$, where $P(t,A)$ also satisfies the following condition:
$$(2)\quad P(t_1+t_2,A)=\int_XP(t_2,A-y)P(t_1,dy)$$
for each $t_1<t_2$ and $t_1+t_2$ in $T$.
\par Then
$$(3)\quad \phi (t_1,x_1,t_2,y) := \int_X P(t_1,x_1,t_2,dx) \exp ({\bf i}y(x)) $$
is the characteristic functional of the transition measure
$P(t_1,x_1,t_2,dx)$ for each $t_1<t_2\in T$ and each $x_1\in X$,
where $X^*_{\bf R}$ notates the topologically dual space of all
continuous $\bf R$-linear real-valued functionals $y$ on $X$, $~y\in
X^*_{\bf R}$. Particularly for $P$ satisfying Conditions $(1)$,
$(2)$ with $t_0=0$ its characteristic functional $\phi $ satisfies
the equalities:
$$(4)\quad \phi
(t_1,x_1,t_2,y)=\psi (t_2-t_1,y)\exp ({\bf i}y(x_1)), $$ where
$$(5)\quad \psi (t,y) := \int_XP(t,dx)\exp ({\bf i}y(x))
\quad \mbox{ and}$$
$$(6)\quad \psi (t_1+t_2,y)=\psi (t_2,y)\psi (t_1,y)$$
for each $t_1<t_2\in T$ and $t_2-t_1\in T$ and $t_1+t_2\in T$
respectively and $y\in X^*_{\bf R}$, $x_1\in X$, since $Z({\cal
A}_{r,C})={\bf C}$.

\par {\bf Remark 2.5. Notation.} If $T$ is a $T_1\cap T_{3.5}$ topological
space, then we denote by $C^0_b(T,H)$ the Banach space of all
continuous bounded functions $f: T\to H$ supplied with the norm:
\par $(1)\quad \| f\|_{C^0}:=\sup_{t\in T} \| f(t)\|_H<\infty $, \\
where $H$ is a Banach space over $\bf R$ which may be also a
two-sided ${\cal A}_{r,C}$-module. If $T$ is compact, then
$C^0_b(T,H)$ is isomorphic with the space $C^0(T,H)$ of all
continuous functions $f: T\to H$.
\par For a set $T$ and a complete locally $\bf R$-convex
space $H$ which may be also a two-sided ${\cal A}_{r,C}$-module
consider the product $\bf R$-convex space $H^T:=\prod_{t\in T}H_t$
in the product topology, where $H_t:=H$ for each $t\in T$.
\par Suppose that ${\sf B}$ is a separating algebra on the space either
$X:=X(T,H)=L^q(T,{\cal B}(T),\lambda ,H)$ or $X:= X(T,H)=C^0_b(T,H)$
or on $X=X(T,H)=H^T$, where $\lambda : {\cal B}(T)\to [0,\infty )$
is a $\sigma $-additive measure on the Borel $\sigma $-algebra
${\cal B}(T)$ on $T$, $1\le q\le \infty $. Consider a random
variable $\xi : \omega \mapsto \xi (t,\omega )$ with values in
$(X,{\sf B})$, where $t\in T$, $\omega \in \Omega $, $(\Omega ,
{\cal R}, P)$ is a measure space with an ${\cal A}_{r,C}$-valued
measure $P$, $~ P: {\cal R}\to {\cal A}_{r,C}$.
\par Events $S_1,\dots ,S_n$ are called independent in total if
$P(\mbox{ }_l\prod_{k=1}^nS_k)=\mbox{ }_l\prod_{k=1}^nP(S_k)$.
Subalgebras ${\cal R}^k\subset {\cal R}$ are said to be independent
if all collections of events $S_k\in {\cal R}^k$ are independent in
total, where $k=1,\dots ,n$, $n\in \bf N$. To each collection of
random variables $\xi _{\gamma }$ on $(\Omega ,{\cal R})$ with
$\gamma \in \Upsilon $ is related the minimal algebra ${\cal
R}_{\Upsilon }\subset \cal R$ for which all $\xi _{\gamma }$ are
measurable, where $\Upsilon $ is a set. Collections $\{ \xi _{\gamma
}: $ $\gamma \in \Upsilon ^l \} $ are called independent if such are
${\cal R}_{\Upsilon ^l}$, where $\Upsilon ^l\subset \Upsilon $ for
each $l=1,\dots ,n,$ $n\in \bf N$. \par For $X=C^0_b(T,H)$ or
$X=H^T$ define $X(T,H;(t_1,\dots ,t_n);(z_1,\dots ,z_n))$ as a
closed submanifold in $X$ of all $f: T\to H$, $f\in X$ such that
$f(t_1)=z_1,\dots ,f(t_n)=z_n$, where $t_1,\dots ,t_n$ are pairwise
distinct points in $T$ and $z_1,\dots ,z_n$ are points in $H$. For
$n=1$ and $t_0\in T$ and $z_1=0$ we denote
$X_0:=X_0(T,H):=X(T,H;t_0;0)$.

\par {\bf Definition 2.6.} Suppose that $H$ is a real Banach
space which also may be a two-sided ${\cal A}_{r,C}$-module.
Consider a random function $w(t,\omega )$ with values in the space
$H$ as a random variable such that:
\par $(1)$ the random variable $\omega (t,\omega )-\omega (u,\omega )$ has
a distribution $\mu ^{F_{t,u}},$ where $\mu $ is an ${\cal
A}_{r,C}$-valued measure on $(X(T,H),{\sf B})$, $\mu ^g(A):=\mu
(g^{-1}(A))$ for $g: X\to H$ such that $g^{-1}({\cal R}_H)\subset
\sf B$ and each $A\in {\cal R}_H$. There by $F_{t,u}$ a $\bf
R$-linear operator $F_{t,u}: X\to H$ is denoted, which is prescribed
by the following formula:
$$F_{t,u}(w):=w(t,\omega )-w(u,\omega )$$
for each $u<t$ in $T$, where ${\cal R}_H$ is a separating algebra of
$H$ such that $F_{t,u}^{-1}({\cal R}_H)\subset \sf B$ for each $u<
t$ in $T$, where $T=[0,b]$ with $0<b<\infty $ or $T=[0,\infty )$,
$\Omega \ne \emptyset $;
\par $(2)$ the vectors $w(t_m,\omega )-w(t_{m-1},\omega )$, ...,
$w(t_1,\omega )-w(0,\omega )$ and $w(0,\omega )$ are mutually
independent for each chosen $0<t_1<...<t_m$ in $T$ and each $m\ge
2$, where $\omega \in \Omega .$ \par Then $ \{ w(t): t\in T \} $ is
called the random function with independent increments, where $w(t)$
is the shortened notation of $w(t,\omega )$.
\par It also may be put \par $(3)$ $w(0,\omega )=0$.

\par {\bf Remark 2.7.} The random function $w(t,\omega )$ satisfying
Conditions 2.6$(1)$-$(3)$ possesses the Markovian property with the
transition measure \par $P(u,x,t,A)=\mu ^{F_{t,u}}(A-x)$.
\par As usually it is put for the expectation $$E_Pf=\int_{\Omega
}f(\omega )P(d\omega )=P^L(f)$$ of a random variable $f: \Omega \to
{\cal A}^h_{r,C}$ whenever this integral exists, where $P=P_{[r]}$
is the ${\cal A}_{r,C}$-valued measure on a measure space $(\Omega
_{[r]}, \mbox{ }_{[r]}{\cal F})$ shortly denoted by $(\Omega ,{\cal
F})$, where $f$ is $({\cal F},{\cal B}({\cal
A}^h_{r,C}))$-measurable, $h\in {\bf N}$, ${\cal B}({\cal
A}^h_{r,C})$ denotes the Borel $\sigma $-algebra on ${\cal
A}^h_{r,C}$. If $P$ is specified, it may be shortly written $E$
instead of $E_P$. If ${\cal G}$ is a sub-$\sigma $-algebra in the
$\sigma $-algebra ${\cal F}$ and if there exists a random variable
$g: \Omega \to {\cal A}^h_{r,C}$ such that $g$ is $({\cal G},{\cal
B}({\cal A}^h_{r,C}))$-measurable and
$$\int_A f(\omega ) P(d\omega ) = \int_Ag(\omega ) P(d\omega
)$$ for each $A\in {\cal G}$, then $g$ is called the conditional
expectation relative to ${\cal G}$ and denoted by $g= E(f|{\cal
G})$.

\par Recall that an operator $J: {\cal A}_{r,C}^n\to {\cal
A}_{r,C}^h$ is called right ${\cal A}_{r,C}$-linear in the weak
sense, if
\par $(1)$ $J(xb+yc)=(Jx)b+(Jy)c$ for each $x$ and $y$ in
${\bf R}^n$ and $b$ and $c$ in ${\cal A}_{r,C}$, where the real
field ${\bf R}$ is canonically embedded into the complexified
Cayley-Dickson algebra ${\cal A}_{r,C}$ as ${\bf R}i_0$, $i_0=1$.
Over the algebra ${\bf H}_{\bf C}={\cal A}_{2,C}$ this gives right
linear operators $J(xb+yc)=(Jx)b+(Jy)c$ for each $x$ and $y$ in
${\cal A}_{2,C}^n$ and $b$ and $c$ in ${\cal A}_{2,C}$, since ${\bf
H}_{\bf C}$ is associative. For short we omit "in the weak sense". A
set of such operators we notate by $L_r({\cal A}_{r,C}^n,{\cal
A}_{r,C}^h)$. Then $$ \| J \| =\sup_{z\ne 0; ~ z\in {\cal
A}^n_{r,C}} \frac{\| Jz \| }{\| z \| },$$ where $z=(z_1,...,z_n)$,
$z_j\in {\cal A}_{r,C}$ for each $j\in \{ 1,...,n \} $, where
$$\| z \| ^2=\sum_{j=1}^n \| z_j \|^2,$$
$ \| a \| ^2 =2 |b|^2+2 |c|^2$ for each $a=b+{\bf i}c$ in ${\cal
A}_{r,C}$ with $b$ and $c$ in ${\cal A}_r$ (see also Remark 2.1
\cite{luocmeasaaca20}).
\par In particular, it is useful to consider the following case: $w=J\xi +p$,
where $\xi $ is a ${\bf R}^{2n}$-valued random variable on a
measurable space $(\Omega _{ [0]}, \mbox{ }_{ [0]}{\cal F})$ and
with a probability measure $P_{[0]}: \mbox{ }_{ [0]}{\cal F}\to
[0,1]$, where $p\in {\cal A}_{r,C}^n$, where ${\bf R}^{2n}$ is
embedded into ${\cal A}_{r,C}^n$ as $i_0{\bf R}^n+i_0{\bf i}{\bf
R}^n$, where $J\in L_r({\cal A}_{r,C}^n,{\cal A}_{r,C}^n)$. This
means that $\xi $ is  $(\mbox{ }_{[0]}{\cal F}, {\cal B}({\bf
R}^{2n}))$-measurable, whilst $w$ is $(\mbox{ }_{[r]}{\cal F}, {\cal
B}({\cal A}_{r,C}^n))$-measurable, where $(\Omega _{ [r]}, \mbox{
}_{ [r]}{\cal F})$ is a measurable space, $P_{[r]}: \mbox{ }_{
[r]}{\cal F}\to {\cal A}_{r,C}$ is a measure. \par Assume that there
is an injection $\theta : (\Omega _{[0]}, \mbox{ }_{ [0]}{\cal
F})\to (\Omega _{[r]},\mbox{ }_{ [r]}{\cal F})$ and $P_{[0]}$ has an
extension ${\sf P}=P_{[0]}^{\theta }$ on $(\Omega _{[r]},\mbox{ }_{
[r]}{\cal F})$ such that $P_{[0]}^{\theta }(\Omega _{[r]}\setminus
\theta (\Omega _{[0]}))=0$, $P_{[0]}^{\theta }(A)=P_{[0]}(\theta
^{-1}(A\cap \theta (\Omega _{[0]}))$ for each $A\in \mbox{ }_{
[r]}{\cal F}$ and $|P_{[r]}|(\Omega _{[r]}\setminus \theta (\Omega
_{[0]}))=0$. Then it may be the case that ${\sf P}$ and $P_{[r]}$
are related by Formulas 2.4$(2)$, 2.4$(3)$ \cite{luocmeasaaca20}
with the help of $U=U_{[r]}=J^2$ and $U_{[0]}=I$ using the ${\cal
A}_{r,C}$-analytic extension. If $f=F(w)$, where $F: {\cal
A}_{r,C}^n\to {\cal A}_{r,C}^h$ is a Borel measurable function then
there exists a Borel measurable function $G: {\bf R}^{2n}\to {\cal
A}_{r,C}^h$ such that $G(\xi )=f$. Therefore if $u: {\cal
A}_{r,C}^h\to {\bf R}$ is a Borel measurable function, then using
Formulas 2.4$(2)$, 2.4$(3)$ \cite{luocmeasaaca20} we put
$$Eu(f)=\int_{\Omega _{[0]}} u(G(\xi (\omega )))P_{[0]}(d\omega ).$$
If $$\int_{A_{[0]}}u(G(\xi (\omega ))P_{[0]}(d\omega
)=\int_{A_{[0]}} g(\theta (\omega ))P_{[0]}(d\omega )$$ for each
$A\in {\cal G}$, where $g: \Omega _{[r]}\to {\bf R}$ is $({\cal G},
{\cal B}({\bf R}))$-measurable, $A_{[0]}= \theta ^{-1}(A\cap \theta
(\Omega _{[0]}))$, $\mbox{}_{[0]}{\cal G}= \theta ^{-1}({\cal G}\cap
\theta (\Omega _{[0]}))$, then $g$ will be called the conditional
expectation of $u(f)$ relative to ${\cal G}$ and denoted by
$E(u(f)|{\cal G})=g$, since ${\sf P}(\Omega _{[r]}\setminus \theta
(\Omega _{[0]}))=0$ and $|P_{[r]}|(\Omega _{[r]}\setminus \theta
(\Omega _{[0]}))=0$, where ${\cal G}$ is a $\sigma $-subalgebra in
$\mbox{ }_{ [r]}{\cal F}$.
\par Henceforth this convention will be used, if some other will not be specified.
\par Let
$L_{r,i}({\cal A}_{r,C}^n,{\cal A}_{r,C}^h)$ denote a family of all
right ${\cal A}_{r,C}$-linear operators $J$ from ${\cal A}_{r,C}^n$
into ${\cal A}_{r,C}^h$ fulfilling the condition \par $(2)$ $J({\cal
A}_r^n)\subset {\cal A}_r^h$.

\par {\bf Theorem 2.8.} {\it
Suppose that either $X=C^0_b(T,H)$ or $X=H^T$, where $H={\cal
A}_{r,C}^n$ with $n\in \bf N$, $2\le r<\infty $, either $T=[0,s]$
with $0<s<\infty $ or $T=[0,\infty )$. Then there exists a family
$\Psi $ of pairwise inequivalent Markovian random functions with
${\cal A}_{r,C}$-valued transition measures of the type $\mu
_{Ut,pt}$ (see Definition 2.4 \cite{luocmeasaaca20}) on $X$ of a
cardinality $card (\Psi )= {\sf c}$, where ${\sf c}=2^{\aleph _0}$,
$0<t\in T$.}
\par {\bf Proof.} Naturally
the algebra ${\cal A}_{r,C}^n=\otimes_{j=1}^n{\cal A}_{r,C}$, if
considered as a linear space over ${\bf R}$, also possesses a
structure of the $\bf R$-linear space isomorphic with ${\bf
R}^{2^{r+1}n}$. Therefore the Borel $\sigma $-algebra ${\cal
B}({\cal A}_{r,C}^n)$ of the algebra ${\cal A}_{r,C}^n$ is
isomorphic with ${\cal B}({\bf R}^{2^{r+1}n})$. So put $P(t,A)=\mu
_{Ut,pt}(A)$ for each $0<t\in T$ and $A\in {\cal B}(H)$, where an
operator $U$ and a vector $p$ are marked satisfying conditions of
Definition 2.4 and 2.3$(\alpha )$ \cite{luocmeasaaca20}.
\par Naturally an embedding of ${\bf R}^n$ into ${\cal A}_{r,C}^n$
exists as $i_0{\bf R}^n$, where $i_0=1$. If $\xi (t)$ is an ${\bf
R}^n$-valued random function, $J$ is a right ${\cal A}_{r,C}$-linear
operator $J: {\cal A}_{r,C}^n\to {\cal A}_{r,C}^n$ satisfying the
condition $J({\cal A}_r^n)\subset {\cal A}_r^n$, $v\in {\cal
A}_{r,C}^n$, then generally $w(t) = J\xi (t)+ vt$ is an ${\cal
A}_{r,C}^n$-valued random function, where $0\le t\in T$, $w(t)$ is a
shortened notation of $w(t,\omega )$.

\par It is well-known, that the
operators $B_j^{\pm 1/2}$ exist (see, for example, Ch. IX, Sect. 13
in \cite{gantmb}), since $B_j$ is positive definite for each $j$. On
the Cayley-Dickson algebra ${\cal A}_r$ the function $\sqrt{a}$
exists (see \S 3.7 and Lemma 5.16 in \cite{ludjms7}). It has an
extension on ${\cal A}_{r,C}$ and its branch such that $\sqrt{a}>0$
for each $a>0$ can be specified by the following. Take an arbitrary
$a=a_0+{\bf i}a_1\in {\cal A}_{r,C}$ with $a_0\in {\cal A}_r$ and
$a_1\in {\cal A}_r$. Put $a_{0,0}=Re (a_0)$, $a_{1,0}=Re (a_1)$,
${a_0}' = a_0-a_{0,0}$, ${a_1}'=a_1-Re (a_1)$. If $a_{0,0}\ne 0$ and
$a_{1,0}\ne 0$, then $a$ can be presented in the form $a=(\alpha
+{\bf i}\beta )(u+{\bf i}v')$ with $\alpha \in \bf R$, $\beta \in
\bf R$, $u\in {\cal A}_r$, $v'\in {\cal A}_r$, $Re (v')=0$.
Therefore in the latter case $\sqrt{a}=\sqrt{\alpha +{\bf i}\beta }
\sqrt{u+{\bf i}v'}$, since ${\bf C}=Z({\cal A}_{r,C})$. If $a$ is
such that $a_{0,0}=0$ and $a_{1,0}\ne 0$ then for $b={\bf i}a$ there
are $b_{0,0}=-a_{1,0}\ne 0$ and $b_{1,0}=0$. On the other hand, for
$a$ with $a_{1,0}=0$ the equation $(\gamma +{\bf i}\delta )^2 =
a_0+{\bf i}{a_1}'$ has a solution with $\gamma $ and $\delta $ in
${\cal A}_r$, since utilizing the standard basis of the complexified
Cayley-Dickson algebra this equation can be written as the quadratic
system in $2^r$ complex variables $\gamma _0+{\bf i}\delta
_0,...,\gamma _{2^r-1}+{\bf i}\delta _{2^r-1}$. The latter system
has a solution  $(\gamma , \delta )$ in ${\cal A}_r^2$, since each
polynomial over $\bf C$ has zeros in $\bf C$ by the principal
algebra theorem. Therefore the initial equation has a solution in
${\cal A}_{r,C}$. Thus the operator
$U^{1/2}=\bigoplus_{j=1}^ma_j^{1/2} B_j^{1/2}$ exists and it
evidently belongs to $L_r({\cal A}_{r,C}^n, {\cal A}_{r,C}^n)$.
\par  Particularly, $J$ can be $J=U^{1/2}$,
while as $\xi (t)$ it is possible to take a Wiener process with the
zero expectation and the unit covariance operator.
\par If $f\in X$, then $T\ni t\mapsto f(t)$ defines a
continuous $\bf R$-linear projection $\pi _t$ from $X$ into $H$.
Therefore, $\pi _{t_n}\times (\pi _{t_{n-1}}\times ...\times \pi
_{t_1})$ provides a continuous $\bf R$-linear projection $\pi _q$
from $X$ into $H^q$ for each $0<t_1<...<t_n\in T$, where $q = \{
t_1,...,t_n \} $. These projections and the Borel $\sigma $-algebras
${\cal B}(H^q)$ on $H^q$ for finite linearly ordered subsets $q$ in
$T$ induce an algebra ${\cal R}(X)$ of $X$. Since $H^T$ is supplied
with the product Tychonoff topology, then a minimal $\sigma
$-algebra ${\cal R}_{\sigma }(H^T)$ generated by ${\cal R}(H^T)$
coincides with the Borel $\sigma $-algebra ${\cal B}(H^T)$. The
topological spaces $T$ and $H$ are separable and relative to the
norm topology on $C^0_b(T,H)$ one gets also ${\cal R}_{\sigma
}(C^0_b(T,H))={\cal B}(C^0_b(T,H))$.
\par By
virtue of Proposition 2.7 \cite{luocmeasaaca20} and Formulas
2.4$(2)$, 2.4$(3)$ \cite{luocmeasaaca20} a characteristic functional
of $P_{U,p}(t,A):=\mu _{Ut,pt}$ fulfills Condition 2.4$(6)$. It is
worth to associate with $P_{U,p}(t,A)$ a spatially homogeneous
transition measure $P_{U,p}(t_1,x_1,t_2,A)$ according to Equation
2.4$(1)$. The representation 2.10$(2)$ \cite{luocmeasaaca20}
implies, that a bijective correspondence exists between $\sigma
$-additive norm-bounded ${\cal A}_{r,C}$-valued measures and their
characteristic functionals, since it is valid for each real-valued
addendum $\mu _{j,k}$ (see, for example, \cite{dalfom,shirb11}) and
$Z({\cal A}_{r,C})={\bf C}$. Moreover, a characteristic functional
of the ordered convolution $(\mu
* \nu )$ of two $\sigma $-additive norm-bounded
${\cal A}_{r,C}$-valued measures $\mu $ and $ \nu $ is the ordered
product $\hat{\mu } \cdot \hat {\nu }$ of their characteristic
functionals $\hat{\mu }$ and $\hat {\nu }$ respectively. Therefore,
Conditions 2.1$(1)$-$(4)$ are satisfied.
\par  Then Formulas 2.1$(5)$, $(7)$, $(8)$ together with the data above
describe an ${\cal A}_{r,C}$-valued Markov cylindrical distribution
$P_{U,p}$ on $X$ (see Corollary 2.6 \cite{luocmeasaaca20} and
Definition 2.1), since $t=t_2-t_1>0$ for each $0<t_1<t_2\in T$. The
space $H$ is Radon by Theorem I.1.2 \cite{dalfom}, since $H$ as the
metric space is separable and complete. From Theorem 2.3 and
Proposition 2.7 \cite{luocmeasaaca20} it follows that $P_{U,p}$ is
uniformly norm-bounded. In view of Theorem 2.15 and Corollary 2.17
\cite{luocmeasaaca20} this cylindrical distribution has an extension
to a norm-bounded measure $P_{U,p}$ on a completion ${\cal R}_P(X)$
of ${\cal R}(X)$, where ${\cal R}_{\sigma }(X)={\cal B}(X)$.
\par Considering different operators $U$ and vectors $p$
and utilizing the Kakutani theorem (see, for example, in
\cite{dalfom}) we infer that there is a family of the cardinality
$\sf c$ of pairwise nonequivalent and orthogonal measures of such
type $P_{U,p}$ on $X$, since each $P$ has the representation
2.10$(2)$ \cite{luocmeasaaca20}.
\par  Let $\Omega =\Omega _{[r]}$ be the set of all elementary events
$$\omega :=\{ f: f\in
X(T,H;(t_0,t_1,\dots ,t_n);(0,x_1,\dots ,x_n)) \} ,
$$
where $\Lambda _{\omega }$ is a finite subset of $\bf N$, $x_i\in
H$, $(t_i: i\in \Lambda _{\omega } )\in \Upsilon _T$ is a subset of
$T\setminus \{ t_0\}$ (see Remark 2.2 and Notation 2.5), where
$t_0=0$, where $ t_i<t_j$ for each $i<j$ in $\Lambda _{\omega }$.
Hence an algebra $\tilde {\sf U}$ exists of cylindrical subsets of
$X_0(T,H)$ induced by the projections $\pi _q: X_0(T,H)\to H^q,$
where $q\in \Upsilon _T$ is a subset in $T\setminus \{ 0 \} $. This
procedure induces the algebra ${\cal R}({\Omega })$ of ${\Omega }$.
So one can consider a Markovian random function corresponding to
$P_{U,p}$ (see Definition 2.6).

\par {\bf Corollary 2.9.} {\it Let $w(t,\omega )$ be a random function given by Theorem 2.8 with the transition measure $\mu
_{Ut,pt}$ for each $t>0$, then \par $(1)$ $E(w(t_2,\omega ) -
w(t_1,\omega ))=(t_2-t_1)p$ and \par $(2)$ $E((w_k(t_2,\omega
)-p_kt_2)(w_h(t_1,\omega )-p_ht_1)) = (t_2-t_1) a_j b_{k-\beta
_{j-1}, h- \beta _{l-1}; j} \delta _{j,l}$ \\
for each $k$ and $h$ in $ \{ 1,...,n \} $, where $0<t_1<t_2\in T$,
$1+\beta _{j-1}\le k \le \beta _j$ and $1+\beta _{l-1}\le h \le
\beta _l$, $j=1,...,m$, $l=1,...,m$, where $E$ means the expectation
relative to $P^L_{U,p}$.}
\par {\bf Proof.} By virtue of Theorem 2.8 the random
function $w(t,\omega )$ has the transition measure \par
$P(t_1,x,t_2,A)=\mu
^{F_{t_2,t_1}}(A-x)=P^L_{(t_2-t_1)U,(t_2-t_1)p}$, where
$x=w(t_1,\omega )$. Therefore Formulas $(1)$ and $(2)$ follow from
Proposition 2.8 and Theorem 2.9 \cite{luocmeasaaca20}.

\par {\bf Definition 2.10.} Let $( \Omega , {\cal F}, P )$ be a measure space
with an ${\cal A}_{r,C}$-valued $\sigma $-additive norm-bounded
measure $P$ on a $\sigma $-algebra ${\cal F}$ of a set $\Omega $
with $P(\Omega )=1$. It is said that there is a filtration $ \{
{\cal F}_t: t\in T \} $, if ${\cal F}_{t_1}\subset {\cal
F}_{t_2}\subset {\cal F}$ for each $t_1<t_2$ in $T$, where ${\cal
F}_t$ is a $\sigma $-algebra for each $t\in T$, where either
$T=[0,s]$ with $0<s<\infty $ or $T=[0,\infty )$. A filtration $ \{
{\cal F}_t: t \in  T \} $ is called normal, if $ \{ B\in {\cal F}:
|P|(B) = 0 \} \subset {\cal F}_0$ and ${\cal F}_t = \bigcap _{T\ni
v>t} {\cal F}_v$ for each $t\in T$. \par Then if for each $t\in T$ a
random variable $u(t): \Omega \to X$ with values in a topological
space $X$ is $({\cal F}_t,{\cal B}(X))$-measurable, then the random
function $ \{ u(t): t \in T \} $ and the filtration $\{ {\cal F}_t:
t\in T \} $ are adapted, where ${\cal B}(X)$ denotes the minimal
$\sigma $-algebra on $X$ containing all open subsets of $X$ (i.e.
the Borel $\sigma $-algebra). Let ${\cal G}$ be a minimal $\sigma
$-algebra on $T\times \Omega $ generated by sets $(v,t]\times A$
with $A\in {\cal F}_v$, also $\{ 0 \} \times A$ with $A\in {\cal
F}_0$. Let also $\mu $ be a $\sigma $-additive measure on $(T\times
\Omega ,{\cal G})$ induced by the measure product $\lambda \times
P$, where $\lambda $ is the Lebesgue measure on $T$. If $u: T\times
\Omega \to X$ is $({\cal G}_{\mu }, {\cal B}(X))$-measurable, then
$u$ is called a predictable random function, where ${\cal G}_{\mu }$
denotes the completion of ${\cal G}$ by $|\mu |$-null sets, where
$|\mu |$ is the variation of $\mu $ (see Definition 2.10 in
\cite{luocmeasaaca20}).
\par The random function given by Corollary 2.9 is called an ${\cal
A}_{r,C}^n$-valued $(U,p)$-random function or shortly $U$-random
function for $p=0$.
\par {\bf Remark 2.11.} Random functions described in the proof
of Theorem 2.8 are ${\cal A}_{r,C}$ generalizations of the classical
Brownian motion processes and of the Wiener processes. \par Let
$w(t)$ be the ${\cal A}_{r,C}^n$-valued $(U,p)$-random function
provided by Theorem 2.8 and Corollary 2.9.  Let a normal filtration
$ \{ {\cal F}_t: t\in T \} $ on $(\Omega , {\cal F}, P)$ be induces
by $w(t)$. Therefore $w(t)$ is $({\cal F}_t,{\cal B}({\cal
A}_{r,C}^n))$-measurable for all $t\in T$; $w(t_1+t_2)-w(t_1)$ is
independent of any $A\in {\cal F}_{t_1}$ for each $t_1$ and
$t_1+t_2$ in $T$ with $t_2>0$. In view of Theorem 2.8 and Corollary
2.9 the conditions ${\sf P}(\Omega \setminus \theta (\Omega
_{[0]}))=0$ and $|P_{[r]}|(\Omega \setminus \theta (\Omega
_{[0]}))=0$ are satisfied, where $\Omega =\Omega _{[r]}$, $ ~ {\cal
F}=\mbox{}_{[r]}{\cal F}$ (see Remark 2.7).
\par Suppose that $ \{ S(t): ~ t\in T \} $ is an $L_r({\cal
A}_{r,C}^n,{\cal A}_{r,C}^h)$ valued random function (that is,
random operator), $S(t)=S(t,\omega )$, $\omega \in \Omega $ (see
also the notation in Remark 2.7). It is called elementary, if a
finite partition $0=t_0<t_1<...<t_k=s$ exists so that \par $(1)$
$S(t) = \sum_{l=0}^{k-1}S_l\cdot ch _{(t_l,t_{l+1}]}$, \\ where
$S_l:\Omega \to L_r({\cal A}_{r,C}^n,{\cal A}_{r,C}^h)$ is $({\cal
F}_l, {\cal B}(L_r({\cal A}_{r,C}^n,{\cal A}_{r,C}^h))$-measurable
for each $l=0,...,k-1$, where $n$ and $h$ are natural numbers, where
$ch _{(t_l,t_{l+1}]}$ denotes the characteristic function of the
segment $(t_l,t_{l+1}] = \{ t\in {\bf R}: t_l<t\le t_{l+1} \} $,
$T=[0,s]$. A stochastic integral relative to $w(t)$ and the
elementary random function $S(t)$ is defined by the formula:
$$(2)\quad \int_0^t S(\tau )dw(\tau ) := \sum_{l=0}^{k-1} S_l
(w(t_{l+1}\wedge t) - w(t_l\wedge t)),$$ where $t\wedge t'= \min
(t,t')$ for each $t$ and $t'$ in $T$. Similarly elementary
$L_{r,i}({\cal A}_{r,C}^n,{\cal A}_{r,C}^h)$ random functions and
their stochastic integrals are defined.  Put \par $(3)$
$<x,y>=x_1\tilde{y}_1+...+x_h\tilde{y}_h$ for each $x$ and $y$ in
${\cal A}_{r,C}^h$, \\ where $y=(y_1,...,y_h)$ with $y_l\in {\cal
A}_{r,C}$ for each $l$, $~\tilde{z}=z_0-z'$ for each $z=z_0+z'$ in
${\cal A}_{r,C}$ with $z_0\in {\bf R}$ and $z'\in {\cal A}_{r,C}$,
$Re (z')=0$. \par Denote by $Q^*$ an adjoint operator of an $\bf
R$-linear operator $Q: {\cal A}_{r,C}^n\to {\cal A}_{r,C}^h$ such
that \par $(4)$ $<Qx,y>=<x,Q^*y>$ for each $x\in {\cal A}_{r,C}^n$
and $y\in {\cal A}_{r,C}^h$. \par Then we put for $Q=A+{\bf i}B$
with $A$ and $B$ in $L_{r,i}({\cal A}_{r,C}^n,{\cal A}_{r,C}^h)$
\par $(5)$ $\| Q\|_2^2= 2 Tr (AA^*)+ 2 Tr (BB^*).$

\par {\bf Lemma 2.12.} {\it Let \par $(i)$ $S(t)$ be an elementary $L_r({\cal A}^n_{r,C},{\cal A}^h_{r,C})$-valued
random variable with $E (\| S(t) \| | {\cal F}_a)<\infty $ ${\sf
P}$-a.e. on $(\Omega ,{\cal F})$ for each $t\in [a,b]$ and let
\par $(ii)$ $w=w_0+{\bf i}w_1$ be an ${\cal A}_{r,C}^n$-valued random function with
$U_0$- and $U_1$- random functions $w_0$ and $w_1$ respectively
having values in ${\cal A}_r^n$ so that $U_0$ and $U_1$ belong to
$L_{r,i}({\cal A}_{r,C}^n,{\cal A}_{r,C}^n)$ and the operator
$U=U_0+{\bf i}U_1$ fulfills Conditions 2.3$(\alpha )$ and of
Definition 2.4 \cite{luocmeasaaca20}, where $w_0$ and $w_1$ are
independent; $0\le a<b<\infty $, $[a,b]\subset T$ (see Definitions
2.10 \cite{luocmeasaaca20}, 2.10 and Remarks 2.7, 2.11 above).
\par Then $E(\int_a^bS(t)dw(t)|{\cal F}_a)=0$ $ ~ {\sf P}$-a.e. on
$(\Omega , {\cal F})$.}
\par {\bf Proof.} This follows from Corollary 2.9$(1)$ and Formulas
2.10$(1)$, 2.10$(2)$, since $0\le (b-a)E (\sum_{l=0}^{k-1} \| S_l \|
| {\cal F}_a)<\infty $ $~ {\sf P}$-a.e. and $E(w(t_2, \omega ) -
w(t_1, \omega ))=0$ for each $t_2>t_1$ in $[a,b]$ for the $U$-random
function $w$.
\par {\bf Lemma 2.13.} {\it Let $S=A+{\bf i}B$, with $A$ and $B$
belonging to $L_{r,i}({\cal A}_{r,C}^n,{\cal A}_{r,C}^h)$, where
$n\in {\bf N}$, $h\in {\bf N}$, $2\le r<\infty $. Then
$$(1)\quad \| S \|_2 ^2= Tr [(A+{\bf i}B)((A^*-{\bf i}B^*)]
+ Tr [(A-{\bf i}B)((A^*+{\bf i}B^*)]<\infty ~ \mbox{ and}$$
$$(2)\quad \| S \| \le \| S \|_2.$$}
\par {\bf Proof.} Since $A$ and $B$
belong to $L_{r,i}({\cal A}_{r,C}^n,{\cal A}_{r,C}^h)$, then \par
$(3)$  $ \| A+{\bf i}B \|_2^2=2Tr (AA^*)+2Tr (BB^*)<\infty $ \\ by
2.11$(5)$, where as usually $Tr (AA^*)$ denotes the trace of the
operator $AA^*$. On the other side,
\par $[(A+{\bf i}B)((A^*-{\bf i}B^*)] + [(A-{\bf i}B)((A^*+{\bf
i}B^*)] =2 (AA^*+BB^*)$. \\
Since $A\in L_{r,i}({\cal A}_{r,C}^n,{\cal A}_{r,C}^h)$, then
$<Ae_k,e_l>\in {\cal A}_r$ for each $k=1,...,n$, $l=1,...,h$, where
$ \{ e_k: k=1,...,m \} $ denotes the standard orthonormal base in
the Euclidean space ${\bf R}^m$, where $m=\max (n,h)$; ${\bf R}^n$
is embedded into ${\cal A}_{r,C}^n$ as $i_0{\bf R}^n$. Therefore we
deduce using Formulas 2.11$(3)$ and 2.11$(4)$ that \par $(3)$ $Tr
(AA^*)=\sum_{l,k}|<e_l,Ae_k>|^2\ge 0$, \\ since $Tr (AA^*)=\sum_l
<AA^*e_l,e_l>=\sum_{l,k}<A^*e_l,e_k> <e_k,A^*e_l>$. \par  This
implies Formula $(1)$. From the Cauchy-Bunyakovskii-Schwarz
inequality, Remark 2.7, Formulas $(1)$ and $(3)$ one gets Inequality
$(2)$.

\par {\bf Theorem 2.14.} {\it
If $S(t)$ is an elementary random function with values in
$L_{r,i}({\cal A}_{r,C}^n,{\cal A}_{r,C}^h)$ and $w(t)$ is an
$U$-random function in ${\cal A}_r^n$ as in Definition 2.10 with
$U\in L_{r,i}({\cal A}_{r,C}^n,{\cal A}_{r,C}^n)$, then
$$(1)\quad E\bigg{[}<\int_a^tS(\tau )dw(\tau ), \int_0^tS(\tau )dw(\tau )>| {\cal F}_a\bigg{]}$$ $$=
E\bigg{[}\int_a^tTr (\{ S(\tau )U^{1/2} \} \{ (U^{1/2})^*S^*(\tau )
\} )d\tau | {\cal F}_a\bigg{]}$$ $~{\sf P}$-a.e. for each $0\le
a<t\in T$.}
\par {\bf Proof.} Since $Ew(t)=0$ and $U: {\cal A}_{r,C}^n\to {\cal
A}_{r,C}^n$, $ ~ U\in L_{r,i}({\cal A}_{r,C}^n,{\cal A}_{r,C}^n)$ by
the conditions of this theorem, then $a_j\in {\cal A}_r\setminus \{
0 \}$ for each $j$ and hence $U^{1/2}: {\cal A}_{r,C}^n\to {\cal
A}_{r,C}^n$ and $U^{1/2}\in L_{r,i}({\cal A}_{r,C}^n,{\cal
A}_{r,C}^n)$, since $U$ satisfies the conditions of Definition 2.4
and 2.3$(\alpha )$ \cite{luocmeasaaca20} (see also Theorem 2.8).
Therefore $w(t,\omega )\in {\cal A}_r^n$ and hence $S(t,\omega )w(t,
\omega )\in {\cal A}_r^h$ for each $t\in T$ and ${\sf P}$-almost all
$\omega \in \Omega $, where $w(t)$ is a shortening of $w(t, \omega
)$, while $S(t)$ is that of $S(t,\omega )$. On the other hand,
$$(2)\quad <x,x>=| x| ^2 = \sum_{j=1}^h
x_j\tilde{x}_j=\sum_{j=1}^h|x_j|^2$$ for each $x\in {\cal A}_r^h$,
where $|z|^2=z\tilde{z}=\sum_{l=0}^{2^r-1}z_l^2$ for each $z$ in the
Cayley-Dickson algebra ${\cal A}_r$, where
$z=z_0i_0+...+z_{2^r-1}i_{2^r-1}$ with $z_l\in {\bf R}$ for each
$l$, $ \{ i_0,...,i_{2^r-1} \} $ is the standard basis of ${\cal
A}_r$.
\par Let $e_l\in {\cal A}_{r,C}^n$ and $f_l\in {\cal A}_{r,C}^h$,
where $e_l=(\delta _{l,k}: k=1,...,n)$ and $f_l=(\delta _{l,k}:
k=1,...,h)$, where $\delta _{l,k}$ is  the Kronecker delta. Then for
an operator $J$ in $L_{r,i}({\cal A}_{r,C}^n,{\cal A}_{r,C}^h)$ and
each $x\in {\cal A}_{r,C}^n$ the representation is valid:
$$(3)\quad Jx=\sum_{k=1}^n\sum_{l=1}^h J_{l,k}x_kf_l,$$ where
$x=x_1e_1+...+x_ne_n$, $x_k\in {\cal A}_{r,C}$ and $~J_{l,k}\in
{\cal A}_r$ for each $k$ and $l$.
\par From the conditions imposed on $U$ (see Definition 2.4
\cite{luocmeasaaca20}) it follows that $U$ and
\par $(4)$ $U^{1/2}=\bigoplus_{l=1}^ma_j^{1/2}B_j^{1/2}$ and
$(U^{1/2})^*=\bigoplus_{l=1}^m\tilde{a}_j^{1/2}B_j^{1/2}$
\\ belong to $L_{r,i}({\cal A}_{r,C}^n,{\cal A}_{r,C}^n)$, since to
the positive definite operator $B_j$ the positive definite matrix
$[B_j]^{1/2}$ with real matrix elements corresponds for each $j$,
also $z^{1/2}\in {\cal A}_r$ for each $z\in {\cal A}_r$.
\par By virtue of Proposition 2.5 and Formulas 2.8$(2)$, $(3)$
\cite{luocmeasaaca20} $\mu _{Ut,0}$ is the ${\cal A}_r$-valued
measure for each $t>0$, since the Cayley-Dickson algebra ${\cal
A}_r$ is power-associative and $\exp_l(z)=\exp (z)$ for each $z\in
{\cal A}_r$.
\par The random function $S(t)w(t)$ is obtained from the
standard Wiener process $\xi $ in ${\bf R}^n$ with the zero
expectation and the unit covariance operator with the help of the
operator $U^{1/2}$ as \par $(5)$ $S(t)w(t)=S(t)U^{1/2}\xi (t)$ \\
according to Theorem 2.8. Therefore, from the Ito isometry theorem
(see it, for example, in \cite{dalfom,gihmskorb}), Formulas
$(2)$-$(5)$ above and Remarks 2.7, 2.11 the statement of this
theorem follows.

\par {\bf Theorem 2.15.} {\it Suppose that
\par $(i)$ $S(t)$ is an elementary
$L_r({\cal A}_{r,C}^n,{\cal A}_{r,C}^h)$ valued random function and
\par $(ii)$ $w=w_0+{\bf i}w_1$ is an ${\cal A}_{r,C}^n$-valued random function
satisfying Condition 2.12$(ii)$, then
$$(1)\quad E\bigg{[} \bigg{\| } \int_a^tS(\tau )dw(\tau )\bigg{\| }^2|{\cal F}_a\bigg{]}
\le \max (\| U_0^{1/2} \|_2^2 , \| U_1^{1/2}\|_2 ^2) ~ E \bigg{[}
\int_a^t\| S(\tau )\|_2^2 d\tau |{\cal F}_a\bigg{]}$$ $~{\sf
P}$-a.e. for each $0\le a<t\in T$.}
\par {\bf Proof.} We consider the following representation
$S(x+{\bf i}y)=(S_{0,0}x)+ (S_{0,1}y)+{\bf i} (S_{1,0}x) +{\bf i}
(S_{1,1}y)$ of $S$ with $S_{l,k}\in L_{r,i}({\cal A}_{r,C}^n,{\cal
A}_{r,C}^h)$ for every $l, k\in \{ 0, 1 \}$ and $z=x+{\bf i}y\in
{\cal A}_{r,C}^n$ with $x$ and $y$ in ${\cal A}_r^n$. For each
$z=x+{\bf i}y\in {\cal A}_{r,C}^n$ we have $|Sz|^2=|(S_{0,0}x)+
(S_{0,1}y)|^2+|(S_{1,0}x) + (S_{1,1}y)|^2$ (see Remark 2.1
\cite{luocmeasaaca20} and Formula 2.14$(2)$ above). On the other
hand, $|v|^2=<v,v>$ for each $v\in {\cal A}_r^h$. For two operators
$G$ and $H$ in $L_{r,i}({\cal A}_{r,C}^n,{\cal A}_{r,C}^h)$ the
inequality is valid $|Tr (GH^*)|^2\le [Tr (GG^*)]\cdot [Tr (HH^*)]$
due to the representation 2.14$(3)$. Applying Theorem 2.14 and Lemma
2.13 (see also Remarks 2.7, 2.11) to $S_{0,0}w_0+
S_{0,1}w_1=(S_{0,0}\oplus S_{0,1})\eta $ and $S_{1,0}w_0 +
S_{1,1}w_1= (S_{1,0}\oplus S_{1,1})\eta $, where $\eta =w_0\oplus
w_1$ and $U=U_0\oplus U_1$, we infer that
$$(2)\quad E\bigg{[} \bigg{\|} \int_a^tS(\tau )dw(\tau )\bigg{\|}^2|{\cal
F}_a\bigg{]}= $$ $$2 E\bigg{[}\int_0^t \bigg{(}\sum_{l,k=0}^1Tr (\{
S_{l,k}(\tau )U^{1/2}_k \} \{ (U^{1/2}_k)^*S^*_{l,k}(\tau ) \}
)\bigg{)}d\tau |{\cal F}_a\bigg{]}$$
$$\le \max (\| U_0^{1/2} \|_2^2 , \| U_1^{1/2}\|_2^2)~  E\bigg{[}
\int_a^t\| S(\tau )\|_2 ^2 d\tau |{\cal F}_a\bigg{]}$$ $~{\sf
P}$-a.e. for each $0\le a<t\in T$, since $|Tr (GH^*)|=|Tr (HG^*)|$
and $|a+b|\le |a|+|b|$ for each $a$ and $b$ in ${\cal A}_r^h$.

\par {\bf Lemma 2.16.} {\it If conditions 2.15$(i)$, 2.12$(ii)$ are satisfied, then
$$\quad (1) \quad {\sf P} \bigg{ \{ } \bigg{ \| } \int_a^b S(t)dw(t)\bigg{ \|} >\beta
\max (\|U_0^{1/2}\|_2, \| U_1^{1/2} \|_2) \bigg{ \} } \le $$
$$\alpha \beta ^{-2} + {\sf P} \bigg{ \{ } \int_a^b\| S(t) \|
_2^2 dt
>\alpha \bigg{ \} } $$ for each $\alpha >0$, $\beta >0$, $[a,b]\subset
T$, $0\le a <b<\infty $. }
\par {\bf Proof.} According to Formula 2.10$(1)$ $S(t)=S(t_l)$ for each $t_l<t\le t_{l+1}$,
where $a=t_0<t_1<...<t_k=b$. Since $S(t)$ is $({\cal F}_{t_l}, {\cal
B}(L_r({\cal A}_{r,C}^n, {\cal A}_{r,C}^h))$-measurable for each
$t\in (t_l, t_{l+1}]$, then $\int_a^{t_{l+1}} \| S(t) \|_2^2 dt$ is
$({\cal F}_{t_l}, {\cal B}([0,\infty ]))$-measurable. We consider a
modified elementary random function $S_{\alpha }(t)$ such that
$S_{\alpha }(t)=S(t)$ for each $t\le t_l$ if $\int_a^{t_{l+1}} \|
S(t) \|_2^2dt \le \alpha $; otherwise $S_{\alpha }(t)=0$ for each
$t\in (t_l,b]$ if $\int_a^{t_l} \| S(t) \|_2^2dt \le \alpha <
\int_a^{t_{l+1}} \| S(t) \|_2^2dt  $ for some $l$. Therefore
$\int_a^t \| S_{\alpha }(t) \|_2^2dt \le \alpha $ for each $t\in
[a,b]$ and hence $$(2)\quad {\sf P} \{ \sup_{t\in[a,b]} \| S_{\alpha
}(t)-S(t)\|_2 >0 \} = {\sf P} \bigg{ \{ } \int_a^b \| S(t) \| _2^2
dt
>\alpha \bigg{ \} }.$$ Then we deduce that
$${\sf P}  \bigg{ \{ } \bigg{ \| } \int_a^b S(t)dw(t)\bigg{ \| } >\beta \max
(\|U_0^{1/2}\|_2, \| U_1^{1/2} \|_2) \bigg{ \} } =$$  $${\sf P}
\bigg{ \{ } \bigg{ \| } \int_a^b S_{\alpha }(t)dw(t)+ \int_a^b
(S(t)-S_{\alpha }(t))dw(t)\bigg{ \| }
>\beta \max (\|U_0^{1/2}\|_2, \| U_1^{1/2} \|_2) \bigg{ \} } \le $$ $$ {\sf P}
\bigg{ \{ } \bigg{ \| }\int_a^b S_{\alpha }(t)dw(t) \bigg{ \| }
>\beta \max (\|U_0^{1/2}\|_2, \| U_1^{1/2} \|_2) \bigg{ \} } + {\sf P}
\bigg{\{ }\bigg{ \|}\int_a^b (S(t)-S_{\alpha }(t))dw(t) \bigg{ \| }
>0 \bigg{ \} }
$$
$$ \le  \frac{E\bigg{[ } \bigg{\| } \int_a^b S_{\alpha
}(t)dw(t)\bigg{\|}^2\bigg{]}}{\beta ^2 \max (\|U_0^{1/2}\|_2^2, \|
U_1^{1/2} \|_2^2)} + {\sf P} \bigg{ \{ } \int_a^b \| S(t) \|_2^2dt
>\alpha \bigg{ \} }
$$ by Chebysh\"ev inequality (see it, for example, in Sect. II.6 \cite{shirb11}), Equality $(2)$ above, Formulas
2.10$(1)$, $(2)$ in \cite{luocmeasaaca20}. By virtue of Theorem 2.15
$$E\bigg { [ }\bigg{ \| } \int_a^b S_{\alpha }(t) dw(t) \bigg{ \| }^2 \bigg{ ] }
\le \max (\|U_0^{1/2}\|_2^2, \| U_1^{1/2} \|_2^2) ~  E\bigg{[  }
\int_a^b  \| S(t) \|_2^2 dt \bigg{ ] },$$ since $E[E(\zeta |{\cal
F}_a)]=E\zeta $ for a random variable $\zeta : \Omega \to [0,\infty
]$ which is $({\cal F}_a, {\cal B}([0,\infty ]))$-measurable (Sect.
II.7 \cite{shirb11}). This implies inequality $(1)$.

\par {\bf Theorem 2.17.} {\it If $w$ is a $U$-random function and $\{ S(t): t\in T \} $
is an $L_{r,i}({\cal A}_{r,C}^n,{\cal A}_{r,C}^h)$-valued
predictable random function satisfying the condition
$$(1)\quad E[\int_a^tTr (\{ S(\tau )U^{1/2} \} \{ (U^{1/2})^*S^*(\tau ) \}
)d\tau ]<\infty $$ for each $0\le a<t$ in $T$, where the operator
$U$ is specified in Definition 2.4 \cite{luocmeasaaca20} such that
$U\in L_{r,i}({\cal A}_{r,C}^n,{\cal A}_{r,C}^n)$, then a sequence
$\{ S_{\kappa }(t): \kappa \in {\bf N} \} $ of elementary random
functions exists with $t\in T$ such that
$$(2)\quad lim_{\kappa \to \infty } E[\int_a^tTr (\{ (S(\tau )-
S_{\kappa }(\tau ))U^{1/2} \} \{ (U^{1/2})^*(S^*(\tau )-S_{\kappa
}(\tau )) \} )d\tau ]=0$$ for each $0\le a<t$ in $T$.}
\par {\bf Proof.} Notice that
$Tr (\{ S(\tau )U^{1/2} \} \{ (U^{1/2})^*S^*(\tau ) \} \ge 0$ for
each $\tau \in T$, since $U\in L_{r,i}({\cal A}_{r,C}^n,{\cal
A}_{r,C}^n)$ implying $a_j\in {\cal A}_r$ and hence $a_j^{1/2}\in
{\cal A}_r$ for each $j$. In view of Formulas 2.14$(1)$, $(3)$ the
random function $S(\tau )U^{1/2}$ having values in $L_{r,i}({\cal
A}_r^n,{\cal A}_r^h)$ has the decomposition into a finite $\bf
R$-linear combination
$$(3)\quad S(t)U^{1/2}=\sum_{l=1}^n\sum_{k=1}^h\sum_{j=0}^{2^r-1}\eta _{l,k;j}
e_l\otimes f_ki_j$$ of real random functions $\eta _{l,k;j}$ using
vectors $e_l$, $f_k$ and the standard basis $ \{ i_0, i_1,...,
i_{2^r-1} \} $ of the Cayley-Dickson algebra ${\cal A}_r$ over ${\bf
R}$. For each real-valued random function the condition
$$(4)\quad E[\int_a^t\eta _{l,k;j}^2 d\tau ]<\infty $$ is fulfilled
for each $0\le a<t$ in $T$, hence a sequence of real-valued random
functions $\eta _{l,k;j;\kappa }$ exists such that
$$(5)\quad lim_{\kappa \to \infty } E[\int_a^t (\eta _{l,k;j}-\eta
_{l,k;j;\kappa })^2 d\tau ]=0$$ for each $t\in T$. Thus Formulas
$(3)$ and $(5)$ imply $(2)$.

\par {\bf Theorem 2.18.} {\it If  $w$ fulfills Condition
2.12$(ii)$ and $S(t)$ is a $L_r({\cal A}_{r,C}^n,{\cal
A}_{r,C}^h)$-valued predictable random function satisfying the
following condition
$$(1)\quad
E \bigg{[} \int_a^b F(S;U_0,U_1)(\tau )d\tau \bigg{]} <\infty $$ for
each $0\le a<b$ in $T$, where
$$(2)\quad F(S;U_0,U_1)(t)= \sum_{l,k=0}^1 Tr (\{ S_{l,k}(t)U^{1/2}_k \} \{
(U^{1/2}_k)^*S^*_{l,k}(t) \} ),$$  then a sequence $\{ S_{\kappa
}(t): \kappa \in {\bf N} \} $ of elementary random functions exists
with $t\in T$ such that
$$(3)\quad lim_{\kappa \to \infty } E\bigg{[}\int_a^bF((S(\tau )-
S_{\kappa }(\tau ));U_0,U_1)(\tau ) d\tau \bigg{]}=0$$ for every
$0\le a<b$ in $T$.}
\par The {\bf proof} is analogous to that of Theorem 2.17 with the help of Formula 2.15$(2)$,
since $E(E(\zeta |{\cal F}_a))=E\zeta $ with $\zeta =
\int_a^bF(S;U_0,U_1)(\tau )d\tau $, $ ~ \zeta \ge 0$ $~ {\sf
P}$-a.e.

\par {\bf Definition 2.19.} It will be said that
a sequence $\{ S_{\kappa }(t): \kappa \in {\bf N} \} $ of elementary
$L_r({\cal A}_{r,C}^n,{\cal A}_{r,C}^h)$-valued random functions
with $t\in T$ is mean absolute square convergent to a predictable
$L_r({\cal A}_{r,C}^n,{\cal A}_{r,C}^h)$-valued random function $\{
S(t): t\in T \} $, where $w$ satisfies Condition 2.12$(ii)$, if
Condition 2.18$(3)$ is fulfilled. The corresponding mean absolute
square limit is induced by Formulas 2.15$(2)$, 2.18$(3)$ and is
denoted by $l.i.m.$. The family of all predictable $L_r({\cal
A}_{r,C}^n,{\cal A}_{r,C}^h)$-valued random functions $\{ S(t): t\in
T \} $ satisfying condition 2.18$(1)$ will be denotes by
$V_{2,1}(U_0,U_1,a,b,n,h)$
\par A stochastic integral of $S\in V_{2,1}(U_0,U_1,a,b,n,h)$ is:
$$(1) \quad \int_0^tS(\tau )dw(\tau ) := l.i.m. _{\kappa \to \infty }
\int_0^tS_{\kappa }(\tau )dw(\tau ),$$ where $(2)$ $w=w_0+{\bf
i}w_1$ is an ${\cal A}_{r,C}^n$-valued random function with $U_0$
and $U_1$ random functions $w_0$ and $w_1$ respectively having
values in ${\cal A}_r^n$, where $0\le a\le t\le b $ in $T$, where
$w$ satisfies Condition 2.12$(ii)$.

\par {\bf Proposition 2.20.} {\it Let the conditions of Theorem
2.18 be satisfied and let $S\in V_{2,1}(U_0,U_1,a,b,n,h)$, $0\le
a<c<b\in T$, then there exists $\int_{\beta }^{\gamma }S(t)dw(t)$
for each $a\le \beta \le \gamma \le b$ and
$$(1)\quad \int_a^b S(t)dw(t)=\int_a^cS(t)dw(t)+\int_c^bS(t)dw(t).$$}
\par {\bf Proof.} In view of Theorem 2.18, Definitions 2.10, 2.19 and Remark 2.11
there exists $\int_{\beta }^{\gamma }S(t)dw(t)$ for each $a\le \beta
\le \gamma \le b$. Formula $(1)$ for elementary random functions
$S_{\kappa }$ for each $\kappa \in {\bf N}$ follows from 2.11$(2)$.
Hence taking $l.i.m._{k\to \infty }$ we infer Equality $(1)$ for
$S\in V_{2,1}(U_0,U_1,a,b,n,h)$ by Theorem 2.18.

\par {\bf Proposition 2.21.} {\it If $S\in V_{2,1}(U_0,U_1,a,b,n,h)$,
$S_{\kappa }\in V_{2,1}(U_0,U_1,a,b,n,h)$ for each $\kappa \in {\bf
N}$, $w$ satisfies Condition 2.12$(ii)$, and $$(1)\quad \lim_{\kappa
\to \infty } E\int_a^bF(S-S_{\kappa };U_0,U_1)(t)dt=0,$$ where $0\le
a<b\in T$, then there exists $$(2)\quad l.i.m._{\kappa \to \infty
}\int_a^bS_{\kappa }(t)dw(t)=\int_a^bS(t)dw(t).$$}
\par {\bf Proof.} In view of Proposition 2.20
stochastic integrals $\int_a^bS(t)dw(t)$ and $\int_a^bS_{\kappa
}(t)dw(t)$ exist for each $\kappa \in {\bf N}$. From Theorem 2.18
and Definition 2.19 it follows that
$$l.i.m._{\kappa \to \infty }\int_a^bS_{\kappa
}(t)dw(t)=\int_a^bS(t)dw(t).$$

\par {\bf Proposition 2.22.} {\it If $S\in V_{2,1}(U_0,U_1,a,b,n,h)$,
and if $w$ satisfies Condition 2.12$(ii)$, where $0\le a <b\in T$,
then
$$(1)\quad E\bigg{[}\int_a^bS(t)dw(t)\bigg{|}{\cal F}_a\bigg{]}=0\quad {\sf P} \mbox{-a.e.
and}$$
$$(2)\quad E\bigg{[} \bigg{\|} \int_a^tS(\tau )dw(\tau )\bigg{\|}^2|{\cal
F}_a\bigg{]}= 2 E\bigg{[}\int_0^t F(S;U_0,U_1)(\tau ) d\tau |{\cal
F}_a\bigg{]}$$
$$\le \max (\| U_0^{1/2} \|_2^2 , \| U_1^{1/2}\|_2^2)~  E\bigg{[}
\int_a^t\| S(\tau )\|_2 ^2 d\tau |{\cal F}_a\bigg{]}$$ $ ~ {\sf
P}$-a.e. for each $0\le a<t\in T$.}
\par {\bf Proof.} From Lemmas 2.12, 2.13, Proposition 2.20 the identity
$(1) $ follows. Then Theorem 2.15 and Proposition 2.20 imply
Inequality $(2)$, since $E(E(\zeta |{\cal F}_a))=E\zeta $ with
$\zeta =\int_a^bF(S;U_0,U_1)(t)dt$ and since \par ${\sf P} \bigg{ \{
} \omega \in \Omega : ~ E\bigg{[} \int_a^b F(S;U_0,U_1)(t)dt\bigg{|}
{\cal F}_a\bigg{]}(\omega )=\infty  \bigg{ \} }=0$; $ ~ \zeta \ge 0$
$~ {\sf P}$-a.e.

\par {\bf Remark 2.23.} Let $ch_{[0,\infty )}(t)=1$ for each $t\ge
0$, and $ch_{[0,\infty )}(t)=0$ for each $t<0$, be a characteristic
function of $[0,\infty )$, $~ [0,\infty )\subset {\bf R}$. Then
$G(\tau )\in V_{2,1}(U_0,U_1,a,b,n,h)$ for each $t\in [a,b]$, if
$S(\tau )\in V_{2,1}(U_0,U_1,a,b,n,h)$, where $G(\tau ) := S(\tau
)ch_{[0,\infty )}(\tau -t)$. It is put $$(1)\quad \eta
(t)=\int_a^tS(\tau )dw(\tau ):=\int_a^bS(\tau )ch_{[0,\infty )}(t-
\tau )dw(\tau )$$ for each $t\in [a,b]$. From Proposition 2.22 it
follows that $\eta (t)$ is defined ${\sf P}$-a.e. By virtue of
Theorem IV.2.1 in \cite{gihmskorb} $\eta (t)$ is the separable
random function up to the stochastic equivalence, since $({\cal
A}_{r,C}^h, |\cdot |)$ is the metric space. Therefore $\eta (t)$
will be considered as the separable random function.

\par {\bf Definition 2.24.} Let $\zeta (t)$, $t\in T$, be a $L_{r,C}^h$-valued random
function adapted to the filtration $\{ {\cal F}_t: t\in T \} $ of
$\sigma $-algebras ${\cal F}_t$ and let $E|\zeta (t)|<\infty $ for
each $t\in T$. If $E(\zeta (t)|{\cal F}_s)=\zeta (s)$ for each $s<t$
in $T$, then the family $ \{ \zeta (t), {\cal F}_t: t\in T \} $ is
called a martingale. If $\zeta (t)\in {\bf R}$ for each $t\in T$ and
$E(\zeta (t)|{\cal F}_s)\ge \zeta (s)$ for each $s<t$ in $T$, then $
\{ \zeta (t), {\cal F}_t: t\in T \} $ is called a sub-martingale.

\par {\bf Lemma 2.25.} {\it Assume that $S(t)\in V_{2,1}(U_0,U_1,a,b,n,h)$ and
$w$ satisfies Condition 2.12$(ii)$, $0\le a<b<\infty $,
$[a,b]\subset T$ and
$$(1)\quad E\bigg{[}\int_a^bF(S;U_0,U_1)(t)dt\bigg{|}{\cal F}_a\bigg{]}<\infty
$$ and $\eta (t)$ is provided by Formula 2.23$(1)$, then $\{ \eta (t), {\cal F}_t: t\in [a,b] \} $ is a martingale
and $\{ | \eta (t) | ^2, {\cal F}_t: t\in [a,b] \} $ is the
sub-martingale.}
\par {\bf Proof.} By virtue of Proposition 2.22 $\eta (t)$ is
$({\cal F}_{t},{\cal B}({\cal A}_{r,C}^h))$-measurable and $E(\eta
(t_2)-\eta (t_1)|{\cal F}_{t_1})=E\bigg{[}\int_{t_1}^{t_2}S(\tau
)d\tau \bigg{|}{\cal F}_{t_1}\bigg{]}=0$ for each $a\le t_1<t_2\le
b$. Hence $\{ \eta (t), {\cal F}_t: t\in [a,b] \} $ is the
martingale.
\par The random function $\eta (t)$ has the decomposition:
$$ (2)\quad \eta (t) = \sum_{k\in \{ 0, 1 \}; ~ j\in \{ 0,1,...,2^r-1 \}; ~ l\in \{ 1,...,h \} }
\eta _{k,j,l} (t) i_j {\bf i}^k e_l$$ with $\eta _{k,j,l} (t) \in
{\bf R}$, for each $k$, $j$, $l$, where $ \{ e_l: l=1,...,h \} $ is
the standard orthonormal basis of the Euclidean space ${\bf R}^h$,
where ${\bf R}^h$ is embedded into ${\cal A}_{r,C}^h$ as $i_0{\bf
R}^h$. Therefore each random function $\eta _{k,j,l} (t)$ is the
martingale. Then $$(3)\quad |\eta (t)|^2=\sum_{k\in \{ 0, 1 \}; ~
j\in \{ 0,1,...,2^r-1 \}; ~ l\in \{ 1,...,h \} } |\eta _{k,j,l}
(t)|^2.$$ By virtue of Theorem 1 and Corollary 2 in Ch. III, Sect. 1
\cite{gihmskorb} and Formula 2.22$(2)$ above $\{ |\eta _{k,j,l}
(t)|^2, {\cal F}_t: t\in [a,b] \} $ is the sub-martingale for each
$k$, $j$, $l$, consequently, $\{ | \eta (t) | ^2, {\cal F}_t: t\in
[a,b] \} $ is the sub-martingale by Formulas $(2)$ and $(3)$.

\par {\bf Lemma 2.26.} {\it Let $S(t)\in V_{2,1}(U_0,U_1,a,b,n,h)$ and
$w$ satisfy Condition 2.12$(ii)$ such that
$$(1)\quad E\bigg{[}\int_a^bF(S;U_0,U_1)(t)dt\bigg{|}{\cal F}_a\bigg{]}<\infty
,$$ then $$(2)\quad  {\sf P} \{ \sup_{t\in [a,b]}
\bigg{|}\int_a^tS(\tau )dw(\tau )\bigg{|}> \beta ~ \max (\|
U_0^{1/2} \|_2, \| U_0^{1/2} \|_2) \bigg{|} {\cal F}_a \} \le $$  $$
\beta ^{-2} E\bigg{[}\int_a^b F(S;U_0,U_1)(t)dt \bigg{|}{\cal
F}_a\bigg{]},$$
$$(3)\quad  {\sf P} \{ \sup_{t\in [a,b]} \bigg{|}\int_a^tS(\tau
)dw(\tau )\bigg{|}> \beta ~ \max (\| U_0^{1/2} \|_2, \| U_0^{1/2}
\|_2) \} \le $$  $$\beta ^{-2} E\bigg{[}\int_a^b F(S;U_0,U_1)(t)dt
\bigg{]}.$$ }
\par {\bf Proof.} From $(2)$ it follows $(3)$. Therefore, it remains
to prove $(2)$. We take an arbitrary partition $a=t_0<t_1<...<t_n=b$
of $[a,b]$. Then we consider $\eta _k := \int_a^{t_k} S(\tau
)dw(\tau )$. In view of Lemma 2.25 $ \{ \eta _l, {\cal F}_{t_l} :
l=1,...,n \} $ is the martingale and $ \{ |\eta _l|^2, {\cal
F}_{t_l} : l=1,...,n \} $ is the sub-martingale.
\par Therefore, from Theorem 5 in Ch. III, Sect. 1 \cite{gihmskorb}
and Formulas 2.10$(1)$, $(2)$ \cite{luocmeasaaca20} we deduce that
$${\sf P} \{ \sup_{0\le l\le n}|\eta _l|>\beta (\max (\|U_0^{1/2}\|_2,
\| U_1^{1/2}\|_2)|{\cal F}_a \} \le \beta ^{-2}E(|\eta _n|^2|{\cal
F}_a)$$ (see also Remark 2.11). Together with Proposition 2.22 above
and the Fubini theorem (II.6.8 \cite{shirb11}) this implies that
$$(4)\quad {\sf P} \{ \sup_{0\le l\le n} \bigg{|}\int_a^{t_l}S(\tau
)dw(\tau )\bigg{|}> \beta ~ \max (\| U_0^{1/2} \|_2, \| U_0^{1/2}
\|_2) |{\cal F}_a \} \le $$
$$\beta ^{-2} E\bigg{[}\int_a^b F(S;U_0,U_1)(t)dt | {\cal
F}_a\bigg{]}.$$ The random function $\int_a^tS(\tau ) dw(\tau )$ is
separable (see Remark 2.23), hence from $(4)$ it follows $(2)$.

\par {\bf Theorem 2.27.} {\it Let $S\in V_{2,1}(U_0,U_1,a,b,n,h)$ be a predictable
$L_r({\cal A}_{r,C}^n,{\cal A}_{r,C}^h)$-valued random function, let
$w$ satisfy Condition 2.12$(ii)$, $[a,b]\subset T$. Then the random
function $\eta(t)=\int_a^tS(\tau )dw(\tau )$ is stochastically
continuous, where $t\in [a,b]$.}
\par {\bf Proof.} If $S_{\kappa }(\tau )$ is an elementary $L_r({\cal A}_{r,C}^n,{\cal A}_{r,C}^h)$-valued
random function, then $\eta_{\kappa }(t)=\int_a^tS_{\kappa }(\tau
)dw(\tau )$ is stochastically continuous by Formula 2.10$(2)$, since
$w(t)$ is stochastically continuous.
\par For each $S\in V_{2,1}(U_0,U_1,a,b,n,h)$ according to
Definition 2.19 and the Fubini theorem
$\int_a^bE(F(S;U_0,U_1))(t)dt<\infty $. By virtue of Theorem 2.18
there exists a sequence $ \{ S_{\kappa }(t): ~ \kappa \in {\bf N} \}
$ of elementary $L_r({\cal A}_{r,C}^n,{\cal A}_{r,C}^h)$-valued
random functions such that 2.18$(3)$ is satisfied. From Lemma 2.26
and the Fubini theorem we infer that
$${\sf P} \{ \sup_{t\in [a,b]} |\int_a^tS(\tau )dw(\tau ) - \int_a^t
S_{\kappa }(\tau )dw(\tau )|> \epsilon \max ( \| U_0^{1/2} \|_2, \|
U_1^{1/2} \|_2 ) \} $$ $$\le \epsilon ^{-2} \int_a^b E(F(S-S_{\kappa
};U_0,U_1))(t)dt.$$ Therefore, there exists a sequence $ \{ \epsilon
_{\kappa }: \kappa \in {\bf N} \} $ with $\lim_{\kappa \to \infty }
\epsilon _{\kappa }=0$ and a sequence $\{ n_k\in {\bf N}: k\in {\bf
N} \} $ such that
$$\sum_{k=1}^{\infty }\epsilon _k^{-2} \int_a^b
E(F(S-S_{n_k};U_0,U_1))(t)dt<\infty ,$$ consequently,
$$\sum_{k=1}^{\infty } {\sf P} \{ \sup_{t\in [a,b]} |\int_a^b S(\tau
)dw(\tau )- \int_a^b S_{n_k}(\tau )dw(\tau )|>\epsilon _k \} <\infty
.$$ In view of the Borel-Cantelli lemma (see, for example, Ch. II,
Sect. 10 \cite{shirb11}) a natural number $k_0\in {\bf N}$ exists
such that
$${\sf P} \{ \sup_{t\in [a,b]} |\int_a^b S(\tau )dw(\tau )- \int_a^b
S_{n_k}(\tau )dw(\tau )|>\epsilon _k \} =1$$ for each $k\ge k_0$.
Hence $\int_a^tS(\tau )dw(\tau )$ is stochastically continuous,
since $\int_a^tS_{n_k}(\tau )dw(\tau )$ is stochastically continuous
for each $k\in {\bf N}$.

\par {\bf Definition 2.28. The generalized Cauchy problem over the
complexified Cayley-Dickson algebra ${\cal A}_{r,C}$.} Let \par
$(1)$ $H: T\times {\cal A}_{r,C}^h\to L_r({\cal A}_{r,C}^n,{\cal
A}_{r,C}^h)$, \par $(2)$ $G: T\times {\cal A}_{r,C}^h\to {\cal
A}_{r,C}^h$ and \par $(3)$ $w=w_0+{\bf i}w_1$ be a random function
in ${\cal A}_{r,C}^n$ satisfying Condition 2.12$(ii)$, where $n$ and
$h$ are natural numbers.
\par A stochastic Cauchy problem over ${\cal
A}_{r,C}$ is: \par $(4)$ $dY(t)=G(t,Y(t))dt+H(t,Y(t))dw(t)$ with $Y(a)=\zeta $, \\
where $Y(t)$ is an ${\cal A}_{r,C}^h$-valued random function, $\zeta
$ is an ${\cal A}_{r,C}^h$-valued random variable which is ${\cal
F}_a$-measurable, $t\in [a,b]\subset T$, $0\le a<b$. The problem
$(4)$ is understood as the integral equation
$$(5)\quad Y(t)=\zeta + \int_a^tG(\tau ,Y(\tau ))d\tau + \int_a^tH(\tau ,Y(\tau ))dw(\tau
) \mbox{, where } t\in [a,b]\subset T.$$ Then $Y(t)$ is called a
solution, if it satisfies the following conditions $(6)$-$(8)$: \par
$(6)$ $Y(t)$ is predictable,
\par $(7)$ $\forall t\in [a,b]$ $~ {\sf P} \{ Y(t): \int_a^t\|
G(\tau ,Y(\tau )) \|d\tau = \infty \} =0$ and
\par $(8)$ ${\sf P} \{ \omega \in \Omega : ~ \exists t\in [a,b], ~ Y(t)
\ne $ \par $
\zeta + \int _a^t G(\tau ,Y(\tau ))d\tau + \int_a^tH(\tau ,Y(\tau ))dw(\tau ) \} =0 $, \\
where $Y(t)$ is a shortening of $Y(t,\omega )$.

\par {\bf Theorem 2.29.} {\it Let $G(t,y)$ and $H(t,y)$ be Borel functions,
let $w$ satisfy Condition 2.12$(ii)$, let $K=const>0$ be such that
\par $(i)$ $\| G(t,x)-G(t,y) \| + \| H(t,x)-H(t,y)\|_2\le K \| x- y
\| $ and
\par $(ii)$ $\| G(t,y) \| ^2 + \| H(t,y) \| _2^2\le K^2 (1+ \| y \|
^2)$ for each $x$ and $y$ in ${\cal A}_{r,C}^h$, $t\in [a,b]=T$,
where $0\le a<b<\infty $,
\par $(iii)$ $E[\|\zeta \|^2]<\infty $. \par Then a solution $Y$ of Equation 2.28$(5)$
exists; and if $Y$ and $Y_1$ are two stochastically continuous
solutions, then \par $(1)$ ${\sf P} \{ \sup_{t\in [a,b]} \| Y(t) -
Y_1(t) \|
>0 \} =0$.}
\par {\bf Proof.} We consider a Banach space $B_{2,\infty }=B_{2,\infty }[a,b]$ consisting of all
predictable random functions $X: [a,b]\times \Omega \to {\cal
A}_{r,C}^h$ such that $X(t)$ is $({\cal F}_t, {\cal B}({\cal
A}_{r,C}^h))$-measurable for each $t\in [a,b]$ and $\sup_{t\in
[a,b]}E [\| X(t) \| ^2 ]<\infty $ with the norm \par $(2)$ $\| X
\|_{B_{2,\infty }}= (\sup_{t\in [a,b]}E [\| X(t) \| ^2 ])^{1/2}$.
\par In view of Proposition 2.21 there exists and operator $Q$ on $B_{2,\infty
}$ such that
$$(3)\quad QX(t)=\zeta + \int_a^tG(\tau ,X(\tau ))d\tau + \int_a^tH(\tau ,X(\tau ))dw(\tau
)$$ for each $t\in [a,b]$, since $G$ and $H$ satisfy Condition
$(ii)$. Then $QX(t)$ is $({\cal F}_t,{\cal B}({\cal
A}_{r,C}^h))$-measurable for each $t\in [a,b]$, since $G$ and $H$
are Borel functions and $X\in B_{2,\infty }$. By virtue of
Proposition 2.22, using the inequality $(\alpha + \beta +\gamma
)^2\le 3(\alpha ^2+\beta ^2+\gamma ^2)$ for each $\alpha $, $\beta $
and $\gamma $ in ${\bf R}$, the Cauchy-Bunyakovskii-Schwarz
inequality, and Condition $(ii)$ of this theorem, we infer that
$$E[\| QX(t) \| ^2]\le 3E[\|\zeta \| ^2 ] + 3(b-a) \int_a^t K^2(1+
\| X(\tau ) \|^2)d\tau  + 3 E\int_a^t K^2(1+ \| X(\tau ) \|^2)d\tau
$$
$$\le
3E[\|\zeta \| ^2 ] + 3K^2 [(b-a)+1] E\int_a^b(1+ \| X(\tau )
\|^2)d\tau $$
$$\le 3E[\|\zeta \| ^2 ] + 3K^2 (b-a) [(b-a)+1] (1+ \| X \|_{B_{2,\infty }}^2).$$
Thus $Q: B_{2,\infty }\to B_{2,\infty }$. Then using the
Cauchy-Bunyakovskii-Schwarz inequality, 2.3$(12)$
\cite{luocmeasaaca20}, Proposition 2.22, Condition $(i)$ of this
theorem, and the inequality $(\alpha +\beta )^2\le 2(\alpha ^2+\beta
^2)$ for each $\alpha $ and $\beta $ in ${\bf R}$, we deduce that
$$E [ \| QX(t)-X_1(t) \|^2]\le 2(b-a) \int_a^t E [ \| G(\tau ,X(\tau
))-G(\tau , X_1(\tau )) \| ^2]d\tau $$ $$+  2 E [ \| \int_a^t \{
H(\tau ,X(\tau )) -H(\tau ,X_1(\tau )) \} dw(\tau ) \|^2] $$
$$\le C_1 \int_a^t E [\| X(\tau ) - X_1(\tau ) \| ]^2 d\tau \le
C_1(t-a)\|X-X_1 \|_{B_{2,\infty }}^2$$ for each $X$ and $X_1$ in
$B_{2,\infty }$, $t\in [a,b]$, where $C_1=2K^2(b-a+1)$. Therefore,
the operator $Q: B_{2,\infty }\to B_{2,\infty }$ is continuous. Then
we infer that $$E[ \| Q^mX(t)-Q^mX_1(t) \|^2]\le C_1 \int_a^tE[ \|
Q^{m-1}X(\tau )-Q^{m-1}X_1(\tau ) \|^2]d\tau $$ $$\le ...\le C_1^m
\int ... \int_{a<t_1<...<t_n<t} E[ \| X(t_m)-X_1(t_m)
\|^2]dt_1...dt_m $$ $$\le C_1^m\| X-X_1\|_{B_{2,\infty
}}^2(b-a)^m/m!$$ for each $X$ and $X_1$ in $B_{2,\infty }$,
$m=1,2,3,...$. Therefore, \par $\|Q^{m+1}X-Q^mX \|^2_{B_{2,\infty
}}\le C_1^m(b-a)^m \| QX-X \|_{B_{2,\infty }}^2/m!$ for each
$m=1,2,3,...$. Hence the series $\sum_{m=1}^{\infty } \|
Q^{m+1}X-Q^mX \|_{B_{2,\infty }}$ converges. Thus the following
limit exists $\lim_{m\to \infty } Q^mX(t)=:Y(t)$ in $B_{2,\infty }$.
From the continuity of $Q$ it follows that $\lim_{m\to \infty
}Q(Q^mX)=QY$, hence $QY=Y$. Thus \par $ \|QY-Y \|_{B_{2,\infty
}}=0$, consequently, ${\sf P} \{ Y(t)=QY(t) \} =1$ for each $t\in
[a,b]$. This means that $Y(t)$ is the solution of 2.28$(5)$. In view
of Theorem 2.27 and Condition $(ii)$ the solution $Y(t)$ is
stochastically continuous up to the stochastic equivalence.
\par Let now $Y$ and $Y_1$ be two stochastically continuous solutions
of Equation 2.28$(5)$. We consider a random function $q_N(t)$ such
that $q_N(t)=1$ if $\| Y(\tau ) \| \le N$ and  $\| Y_1(\tau ) \| \le
N$ for each $\tau \in [a,t]$, $q_N(t)=0$ in the contrary case, where
$t\in [a,b]$, $N>0$. Therefore $q_N(t)q_N(\tau )=q_N(t)$ for each
$\tau <t$ in $[a,b]$, consequently,
$$q_N(t)[Y(t)-Y_1(t)]=q_N(t)[\int_a^tq_N(\tau )
[G(\tau ,Y(\tau ))-G(\tau ,Y_1(\tau ))]d\tau $$ $$+ \int_a^tq_N(\tau
) [H(\tau ,Y(\tau ))-H(\tau ,Y_1(\tau ))]dw(\tau )].$$ On the other
hand, \par $q_N(\tau )[\| G(\tau ,Y(\tau ))-G(\tau ,Y_1(\tau ))\| +
\| H(\tau ,Y(\tau ))-H(\tau ,Y_1(\tau ))\| ] $\par $ \le K q_N(\tau
) \| Y(\tau )-Y_1(\tau ) \| \le 2KN$ \\ by Condition $(i)$. This
implies that $E [q_N(t) \| Y(t)-Y_1(t) \| ^2]<\infty $. Then using
the Fubini theorem, 2.3$(12)$ \cite{luocmeasaaca20}, Proposition
2.22, Lemma 2.26, we deduce that
$$E[q_N(t) \| Y(t)-Y_1(t) \| ^2]\le 2E [q_N(t) \|
\int_a^t q_N(\tau )[G(\tau ,Y(\tau ))-G(\tau ,Y_1(\tau ))d\tau
\|^2]+$$
$$2 E[q_N(t)\| \int_a^tq_N(\tau )[H(\tau ,Y(\tau ))-H(\tau ,Y_1(\tau
))]dw(\tau )\|^2]$$
$$\le 2(b-a)\int_a^tE[q_N(\tau )\| G(\tau ,Y(\tau ))-G(\tau ,Y_1(\tau ))\| ^2]
d\tau $$ $$+ 4 \int_a^tE[q_N(\tau )F(H(\tau ,Y(\tau ))-H(\tau
,Y_1(\tau ));U_0,U_1)]d\tau $$ $$\le 2 K^2 [b-a + \max ( \|
U_0^{1/2} \| _2^2, \| U_1^{1/2} \|_2^2)] \int_a^t E [q_N(\tau ) \|
Y(\tau )-Y_1(\tau ) \|^2]d\tau .$$ Thus a constant $C_2>0$ exists
such that
$$E[q_N(t) \| Y(t)-Y_1(t) \| ^2]\le C_2 \int_a^t E [q_N(\tau )\|
Y(\tau )-Y_1(\tau ) \| ^2]d \tau .$$ The Gronwall inequality
\cite{gihmskorb,gulcastb} implies that $E[q_N(t) \| Y(t)-Y_1(t) \|
^2]=0$, consequently, \par ${\sf P} \{ Y(t)\ne Y_1(t) \} \le {\sf P}
\{ \sup_{t \in [a,b]} \| Y(t) \|>N \} + {\sf P} \{ \sup_{t \in
[a,b]} \| Y_1(t) \|>N \} $. The random functions $Y(t)$ and $Y_1(t)$
are stochastically continuous, hence stochastically bounded,
consequently, $\lim_{N\to \infty } {\sf P} \{ \sup_{t \in [a,b]} \|
Y(t) \|>N \}=0$ and $\lim_{N\to \infty } {\sf P} \{ \sup_{t \in
[a,b]} \| Y_1(t) \|>N \} =0$. Therefore, the random functions $Y(t)$
and $Y_1(t)$ are stochastically equivalent. Thus ${\sf P} \{
\sup_{t\in [a,b]} \| Y(t)-Y_1(t) \| >0 \} =0$.

\par {\bf Corollary 2.30.} {\it Let operators $G$ and $H$ be $G\in L_r({\cal
A}_{r,C}^h,{\cal A}_{r,C}^h)$ and $H\in L_r({\cal A}_{r,C}^n,{\cal
A}_{r,C}^h)$ such that $G$ be a generator of a semigroup $ \{ S(t):
t \in [0, \infty ) \} $. Let also $w(t)$ be a random function
fulfilling Condition 2.12$(ii)$. Then the Cauchy problem $$(1)\quad
Y(t)=\zeta + \int_0^tGY(\tau )d\tau + \int_0^tHdw(\tau ),$$  where
$t\in T$, $E[\| \zeta \|^2]<\infty $, has a solution $$(2)\quad Y(t)
= S(t)\zeta +\int_0^t S(t-\tau )Hdw(\tau )$$ for each $0\le t\in
T$.}

\par {\bf Proof.} The condition $G\in L_r({\cal
A}_{r,C}^h,{\cal A}_{r,C}^h)$ implies that \par $\| G \| =
\sup_{x\in {\cal A}_{r,C}^h, \| x \| =1} \| Gx \| <\infty $, \\
where $ \| x \|^2=\| x_1 \| ^2+...+\|x_h\|^2$, $x=(x_1,...,x_h)\in
{\cal A}_{r,C}^h$, $x_k\in {\cal A}_{r,C}$ for each $k$. As a
realization of the semigroup $S(t)$ it is possible to take $\{
S(t)=\exp_l(Gt): ~ t\ge 0 \} $, since $G$ is a bounded operator and
$\| \exp_l(Gt) \| \le \exp ( \| G \| t)$  for each $t\ge 0$ by
Formulas 2.1$(9)$ and 2.3$(12)$ \cite{luocmeasaaca20}. Therefore
from Theorem 2.29 the assertion of this corollary follows.

\par {\bf Theorem 2.31.} {\it Let $G$, $H$ and $w$ satisfy conditions of Theorem 2.29, let
$Y_{t,z}(t)$ be an ${\cal A}_{r,C}^h$-valued random function
satisfying the following equation:
$$(1)\quad Y_{t,z}(t_1)=z+\int_t^{t_1}G(\tau ,Y_{t,z}(\tau ))d \tau
+ \int_t^{t_1}H(\tau ,Y_{t,z}(\tau ))dw(\tau ) ,$$ where $z\in {\cal
A}_{r,C}^h$, $t<t_1$ in $[a,b]\subset T$, $0\le a<b<\infty $. Then
the random function $Y$ satisfying Equation 2.28$(5)$ is Markovian
with the transition measure \par $(2)$ $P(t,z,t_1,A)=P\{
Y_{t,z}(t_1)\in A \} $ for each $A\in {\cal B}({\cal A}_{r,C}^h)$.}
\par {\bf Proof.} The random function $Y(t)$ is $({\cal F}_t,{\cal B}({\cal
A}_{r,C}^h))$-measurable for each $t\in [a,b]$. On the other hand,
$Y_{t,z}(t_1)$ is induced by the random function $w(t_1)-w(t)$ for
each $t_1\in (t,b]$, where $w(t_1)-w(t)$ is independent of ${\cal
F}_t$. Therefore $Y_{t,z}(t_1)$ is independent of $Y(t)$ and each
$A\in {\cal F}_t$. By virtue of Theorem 2.29 $Y(t_1)$ is unique (up
to stochastic equivalence) solution of the equation
$$(3)\quad Y(t_1)=Y(t)+\int_t^{t_1}G(\tau ,Y(\tau ))d \tau
+ \int_t^{t_1}H(\tau ,Y(\tau ))dw(\tau ) $$ and $Y_{t,Y(t)}(t_1)$
also is its solution, consequently, ${\sf P} \{
Y(t_1)=Y_{t,Y(t)}(t_1) \} =1$. \par Let $f\in C^0_b({\cal
A}_{r,C}^h,{\cal A}_{r,C})$, where $C^0_b({\cal A}_{r,C}^h,{\cal
A}_{r,C})$ denotes the family of all bounded continuous functions
from ${\cal A}_{r,C}^h$ into ${\cal A}_{r,C}$. Let $g\in R_b(\Omega
,{\cal A}_{r,C})$, where $R_b(\Omega ,{\cal A}_{r,C})$ denotes the
family of all random variables $g:\Omega \to {\cal A}_{r,C}$ such
that there exists $C_g=const
>0$ for which \par ${\sf P} \{ \| g \| <C_g \} =1$, where $C_g$ may depend on $g$. We put
\par $(4)$ $q(z,\omega )=f(Y_{t,z}(t_1,\omega ))$, hence
$f(Y(t_1,\omega ))=q(Y(t_1),\omega )$, \\ where $Y(t)$ is a
shortening of $Y(t,\omega )$ as above, $\omega \in \Omega $. Assume
at first that $q$ has the following decomposition:
\par $(5)$ $q(z,\omega )=\sum_{k=1}^mq_k(z)u_k(\omega )$, \\ where
$q_k: {\cal A}_{r,C}^h\to {\cal A}_{r,C}$, $u_k: \Omega \to {\cal
A}_{r,C}$, $m\in {\bf N}$. This implies that $u_k(\omega )$ is
independent of ${\cal F}_t$ for each $k$. Therefore we deduce that
\par $E[g\sum_{k=1}^mq_k(Y(t))u_k(\omega
)]=\sum_{k=1}^mE[gq_k(Y(t))]Eu_k(\omega )$
\par $=E[\sum_{k=1}^mgq_k(Y(t))]Eu_k(\omega )$ \\ and
\par $E[\sum_{k=1}^mq_k(Y(t))u_k(\omega
)|Y(t)]=\sum_{k=1}^mq_k(Y(t))Eu_k(\omega )$, \\consequently,
\par $(6)$ $Egf(Y(t_1))=EgE[f(Y(t_1))|Y(t)]$ \\ for $q$ of the form
$(5)$. This implies that
\par $(7)$ $E[f(Y(t_1))|{\cal F}_t] =v(Y(t))$,  \\ where
$v(z)=Ef(Y_{t,z}(t_1))$.
\par Then $E[ \| g q(Y(\tau ),\omega )\| ^2 ] \le C_g^2\| f \|_C^2$ for each
$\tau \in [a,b]$ by 2.3$(12)$ \cite{luocmeasaaca20}, since $g$ and
$f$ are bounded, where $ \| f \|_C :=\sup_{z\in {\cal A}_{r,C}^h} \|
f(z) \| <\infty $. Therefore for each $\epsilon >0$ there exists
$f_{(\epsilon )}\in C^0_b({\cal A}_{r,C}^h,{\cal A}_{r,C})$ for
which $q_{(\epsilon )}(z,\omega )=f_{(\epsilon )}(Y_{t,z}(t_1,\omega
))$ has the decomposition of type $(5)$ and such that $E[ \|
q_{(\epsilon )}(Y(t),\omega )-q(Y(t),\omega ) \| ^2]<\epsilon
/C_g^2$. Taking $\epsilon \downarrow 0$ one gets that Formulas $(6)$
and $(7)$ are accomplished for each $f\in C^0_b({\cal
A}_{r,C}^h,{\cal A}_{r,C})$. Therefore $P \{ Y(t_1)\in A|{\cal F}_t
\} = P \{ Y(t_1)\in A| Y(t) \} $ for each $A\in {\cal B}({\cal
A}_{r,C}^h)$, $t<t_1$ in $[a,b]$, since the families $R_b(\Omega
,{\cal A}_{r,C})$ and $C^0_b({\cal A}_{r,C}^h,{\cal A}_{r,C})$ of
all such $g$ and $f$ separate points in ${\cal A}_{r,C}^h$. This
implies that $P \{ Y(t_1)\in A| {\cal F}_t \} = P_{t,Y(t)}(t_1,A)$
for each $A\in {\cal B}({\cal A}_{r,C}^h)$, where $P_{t,z}(t_1,A)=P
\{ Y_{t,z}(t_1)\in A \} $.
\par {\bf Conclusion 2.32.} The obtained in this paper results
open new opportunities for integration of PDEs of order higher than
two. Indeed, a solution of a stochastic PDE with real or complex
coefficients of order higher than two can be decomposed into a
solution of a sequence of PDEs of order one or two with ${\cal
A}_{r,C}$ coefficients \cite{ludrend14,ludcmft12}. They can be used
for further studies of random functions and integration of
stochastic differential equations over octonions and the
complexified Cayley-Dickson algebra ${\cal A}_{r,C}$. It is worth to
mention that equations of the type 2.28$(5)$ are related with
generalized diffusion PDEs of the second order. For example, this
approach can be applied to PDEs describing unsteady heat conduction
in solids \cite{ozsa19}, fourth order Schr\"odinger or Klein-Gordon
type PDEs.
\par Another application of obtained results is for the
implementation of the plan described in \cite{luocmeasaaca20}. It is
related with investigations of analogs of Feynman integrals over the
complexified Cayley-Dickson algebra ${\cal A}_{r,C}$ for solutions
of PDEs of order higher than two.

\end{document}